\documentclass[12pt,reqno,a4paper]{amsart}
\usepackage{amsmath,amssymb,amsfonts}
\usepackage[hypertex]{hyperref}
\usepackage[mathscr]{eucal}

\newtheorem{proposition}{Proposition}
\newtheorem{theorem}{Theorem}
\newtheorem{corollary}{Corollary}
\newtheorem{lemma}{Lemma}
\newtheorem*{claim}{Claim}

\theoremstyle{definition}
\newtheorem{remark}{Remark}[section]
\newtheorem{example}{Example}[section]
\newtheorem{definition}{Definition}

\newcommand{\cal}{\EuScript}

\renewcommand{\leq}{\leqslant}

\DeclareMathOperator{\Ber}{Ber}

 \DeclareMathOperator{\Tr}{Tr}

\DeclareMathOperator{\grad}{grad}

\renewcommand{\L}{{\cal L}}

\DeclareMathOperator{\Vol}{Vol}
\DeclareMathOperator{\Conx}{Conx}

\newcommand{\der}[2]{{\frac{\partial {#1}}{\partial {#2}}}}
\newcommand{\lder}[2]{{\partial {#1}/\partial {#2}}}
\newcommand{\dder}[3]{{\frac{\partial^2 {#1}}{\partial {#2}\partial {#3}}}}

\newcommand{\RR}{\mathbb R}

\newcommand{\p}{\partial}

\newcommand{\f}{\mathbf{f}}

\def\a{\alpha}
\def\A{{\bf A}}
\def\e{{\bf e}}
\def\f{{\bf f}}

\def\vare{\varepsilon}
\def\s{\sigma}
\def\ss {{\bf s}}

\newcommand{\D}[1] {|D({#1})|}
\def\DD {{\cal D}}
\def\G {\Gamma}
\def\o{\omega}
\def\F {{\cal F}}
\def\X {{\bf  X}}
\def\Y {{\bf  Y}}
\def\SS {{\bf S}}
 \def\sing {\text{sing}}

\def\grarrow {{\g\,{\buildrel \X\over \longrightarrow}\, \g^\pr}}

\renewcommand{\r}{{\rho}}
\renewcommand{\l}{{\lambda}}

\def\d{\delta}

\def\pr{\prime}
\def\uppernabla {\,^{^S}\!\nabla}
\newcommand{\bs}{{\boldsymbol{s}}}
\newcommand{\FF}{{\boldsymbol{F}}}

\newcommand{\rh}{{\boldsymbol{\rho}}}
\newcommand{\g}{{\boldsymbol{\gamma}}}

\renewcommand{\div}{\mathop{\mathrm{div}}}

\title[Geometry of differential operators of second order]{Geometry of differential operators of second order, the algebra of densities, and groupoids}
\author{H.~M.~Khudaverdian}
\address{School of Mathematics,  University of Manchester, Oxford Road,  Manchester,   M13 9PL,  UK}
\email{khudian@manchester.ac.uk}
\author{Th.~Th. Voronov}
\address{School of Mathematics,  University of Manchester, Oxford Road,  Manchester,   M13 9PL,  UK}
\email{theodore.voronov@manchester.ac.uk}

\date{1 (14) October  2012}
\begin{document}

\begin{abstract}
In our previous works, we introduced, for each (super)manifold, a commutative algebra of densities. It is endowed with a natural invariant scalar product. In this paper, we study geometry of differential operators of second order on this algebra. In the more
conventional language  they correspond to certain operator pencils. We consider the self-adjoint operators  and analyze the  operator pencils that pass through a given operator acting on densities of a particular weight. There are
`singular values' for  pencil parameters. They are related with  interesting geometric  picture.  In particular, we obtain operators  that depend on certain equivalence classes of connections (instead of connections as such). We study the corresponding  groupoids. From this point of view we analyze two examples: the  canonical Laplacian on an odd symplectic supermanifold appearing in   Batalin--Vilkovisky geometry and the Sturm--Liouville operator on the line, related with classical constructions of projective geometry. We also consider the canonical second order semi-density arising on odd symplectic supermanifolds, which has some  similarity with  mean curvature of surfaces in Riemannian geometry.
\end{abstract}

\maketitle


\section{Introduction}
Differential operators of second order  appear in various  problems of mathematical physics.
        The condition that an operator respects the geometric  structure
of the problem under consideration usually fixes this operator almost
uniquely or at least provides a great deal of information about it.
      For example, the standard Laplacian ${\p^2\over \p x^2}+
{\p^2\over \p y^2}+{\p^2\over \p z^2}$
     in Euclidean space $\mathbb E^3$ is defined uniquely
   (up to a constant) by the condition that it is invariant with respect
       to the isometries of $\mathbb E^3$.
     Consider an arbitrary second order  operator
              \begin{equation}\label{secorderoperator}
       \Delta={1\over 2}\left(S^{ab}\p_a\p_b+T^a\p_a+R\right)
                    \end{equation}
      acting on functions on  a manifold $M$. It  defines on $M$ a symmetric
contravariant tensor $S^{ab}$ (its principal symbol).
     For example, for a Riemannian manifold one can take an operator with the principal symbol
  $S^{ab}=g^{ab}$,
where $g^{ab}$ is the metric tensor (with upper indices).  One can fix the scalar $R=0$ in \eqref{secorderoperator}
by the natural condition $\Delta 1=0$.  What about
the first order term $T^a\p_a$ in the operator \eqref{secorderoperator}?
 One can see
that the Riemannian structure can be used to fix  this term as well.
Indeed, consider on $M$ the divergence operator
\begin{equation}\label{firstequation}
   {\div}_\rh \X=\frac{1}{\rho(x)}\,\p_a
                 \bigl(\rho(x)X^a\bigr)\,,
\end{equation}
where $\rh=\rho(x)|D(x)|$ is an arbitrary non-vanishing volume form, and choose $\rh=\rh_{g}=\sqrt {\det g}\D{x}$.
On a Riemannian manifold this volume form is defined uniquely
(up to a constant factor) by the covariance condition.
Thus we arrive at the  second order operator   $\Delta_g$, where, for an arbitrary function $f$,
\begin{multline}\label{laplacebeltrami}
    \Delta_g f= \frac{1}{2}\,{\div}_{_{\rh_{_g}}} {\grad\,}f=
      {1\over 2}{1\over \rho(x)}{\p \over \p  x^a}
    \left(\rho_{_g}(x)g^{ab}{\p f(x)\over \p  x^b}\right)=\\
            {1\over 2}\,
                       \Bigl(
               \p_a\left(g^{ab}\p_b f\right)+
      \p_a\log \rho_{_g}(x) g^{ab}\p_b f(x)
                      \Bigr)
                      ={1\over 2}\left(g^{ab}\p_a\p_bf+
               \left(\p_a g^{ab} +\p_a\log \sqrt {\det g} g^{ab}\right)\p_b f(x)
                      \right)\,.
\end{multline}
We see that a Riemannian structure on a manifold naturally defines a unique
(up to a constant factor)
second order operator on functions $\Delta_g$, called the Laplace--Beltrami operator.
 For this operator, the terms with the first  derivatives
contain  a {\it connection $\nabla$ on volume forms}.
This connection is defined by setting, for an arbitrary volume form,
$\rh=\rho(x)\D{x}$, $\nabla_\X\rh=\p_\X \left({\rh\over \rh_g}\right)\rh_g$, hence
\begin{multline}\label{exampleofflatconnection}
    \nabla_\X\rh=X^a\left(\p_a+\gamma_a\right)\rho(x) |D(x)|=
                X^a\p_a\left({\rh\over \rh_g}\right)\rh_g=
                X^a\p_a\left({\rho(x)\over \sqrt {\det g}}\right)\sqrt {\det g}\D{x}=\\
     X^a\left(\p_a\rho -\p_a\log \left(\sqrt {\det g}\right)\right) \D{x}\,.
\end{multline}
Here the  connection coefficients are given by $\gamma_a=-\p_a(\rho_g(x))=- \p_a\log \sqrt {\det g} $.

Consider another example. Let $S^{ab}(x)$ be an arbitrary symmetric tensor field (not necessarily non-degenerate)
on a manifold $M$ equipped with an affine   connection
$\nabla\colon  \nabla_a\p_b=\G_{ab}^c\p_c$.  An affine connection defines  the second order operator
 $S^{ab}\nabla_a\nabla_b=S^{ab}\p_a\p_b+\dots$.   The principal symbol of this operator is
 the tensor field $S^{ab}(x)$.
An affine connection induces a connection   on
volume forms by the relation $\gamma_a=-\Gamma^b_{ab}$.
In the case of a Riemannian manifold, the tensor $S^{ab}$ can be fixed by the Riemannian metric,
$S^{ab}=g^{ab}$, and the Levi-Civita Theorem provides a unique
symmetric affine connection which preserves the Riemannian structure.
Then we arrive again at the Laplace--Beltrami operator \eqref{laplacebeltrami}.

It is often important to consider differential operators on densities
of an arbitrary weight  $\lambda$.
For example, a density of weight $\l=0$ is an ordinary function,
a volume form is a density of weight $\lambda=1$.
Wave function in Quantum Mechanics can be naturally considered as a
half-density, i.e., a density of  weight $\lambda={1\over 2}$.                                                 

The study of differential operators on densities of arbitrary
weights is a source of  beautiful geometric constructions.
(See, e.g., the works  \cite{{ManZag}, {DuvOvs}, {Lecomte}} and the book \cite{OvsTab}.)
   For example, consider ${\rm Diff} (M)$-modules which appear when we study operators on densities. Here ${\rm Diff} (M)$ is the group of diffeomorphisms.
Let $\DD_{\l}(M)$ be the space of  second order differential operators
  acting on densities of weight $\l$.
  This space has a natural structure of a  Diff $(M)$-module.
  In~\cite{DuvOvs}, Duval and Ovsienko classified these modules for
  all values of $\l$.
  In particular, they wrote down explicit expressions for ${\rm Diff} (M)$-isomorphisms (intertwining operators)
  $\varphi_{\l\mu}$ between arbitrary modules $\DD_\l(M)$ and $\DD_\mu(M)$
  for $\l,\mu\not=0,{1\over 2},1$.
  These isomorphisms
  have the following form.
    If an operator $\Delta_\l\in \DD_\l(M)$ is given in local coordinates by the expression
       $\Delta_\l=A^{ij}(x)\p_i\p_j+A^i(x)\p_i+A(x)$, then
  its image $\varphi_{\l\mu}(\Delta_\l)=\Delta_\mu\in \DD_\mu(M)$  is given in
  the same local coordinates by the expression
  $\Delta_\mu=B^{ij}(x)\p_i\p_j+B^i(x)\p_i+B(x)$,
    where
    \begin{equation}\label{isomorphism}
    \begin{cases}
    B^{ij}&=A^{ij}\,,\cr
    B^i   &=\frac{2\mu-1}{2\l-1}\,A^i+\frac{2(\l-\mu)}{2\l-1}\,\p_j A^{ji}\,,\cr
    B &=\frac{\mu(\mu-1)}{\l(\l-1)}\,A+\frac{\mu(\l-\mu)}{(2\l-1)(\l-1)}\,
    \left(\p_jA^j-\p_i\p_j A^{ij}\right)\,.\cr
       \end{cases}
    \end{equation}
At the exceptional cases $\l,\mu=0,{1\over 2},1$, non-isomorphic modules occur.

In work~\cite{KhVor2}, we suggested a new approach to studying differential operators on densities.  We defined a commutative algebra $\F(M)$ consisting of all (formal sums of) densities of arbitrary real weights on a manifold $M$. Instead of studying operators $\F_{\l}(M)\to \F_{\mu}(M)$ acting on densities of some particular weight $\l$ and mapping them to densities of some other weight $\mu=\l+\delta$ ($\delta$ being the weight of the operator), we suggest to study differential operators on the algebra $\F(M)$ (that is, differential operators on a commutative algebra).
This algebra possesses a
canonical invariant scalar product.
Due to the existence of the canonical scalar product
it is possible to consider the notions of
self-adjoint and antiself-adjoint differential operators
on this algebra.  Operators on the algebra $\F(M)$ can be identified with pencils of operators
acting between   spaces of densities of various weights. A general differential operator of order $\leq N$ on the algebra $\F(M)$ can be written as a power expansion
\begin{equation*}
    A=A^{(N)}+A^{(N-1)}\,\hat \l +A^{(N-2)}\,\hat \l^2+\ldots+ A^{(1)}\,\hat\l^{N-1}+A^{(0)}\,\hat\l^N\,,
\end{equation*}
where $\hat\l$ is the \emph{weight operator}, whose eigenspaces are the spaces $\F_{\l}(M)$, with the eigenvalues $\l$,  for all $\l\in\RR$. The operator  $\hat\l$ is a derivation of the algebra $\F(M)$, and the coefficients  $A^{(n)}$ in the expansion are usual differential operators    of orders $\leq n$ acting on densities. The corresponding operator pencil is
\begin{equation*}
    A_{\l}=A^{(N)}+A^{(N-1)}\,  \l +A^{(N-2)}\,  \l^2+\ldots+ A^{(1)}\, \l^{N-1}+A^{(0)}\, \l^N\,,
\end{equation*}
where $\l$ is the parameter of the pencil. The evaluation at a particular $\l$, $\l=\l_0$,   gives an operator $A_{\l_0}$ acting on densities of weight $\l_0$. Self-adjoint operators on the algebra $\F(M)$ correspond to pencils satisfying an extra condition~\eqref{selfadjointforpencil} and we refer to them as to `self-adjoint pencils'. In particular,
self-adjoint operators of the second order on $\F(M)$ correspond to certain
``canonical'' operator pencils, which, for a given principal symbol,   are associated with   connections on the bundle of volume forms on $M$ (more precisely,   ``upper connections'').

This approach was put forward and developed in~\cite{KhVor2} for  studying  and classifying second order odd operators on odd symplectic manifolds (arising in  the Batalin--Vilkovisky formalism) and on odd Poisson manifolds.

The  pencils of second order operators that we discovered in~\cite{KhVor2}, which we call ``canonical'',  possess the following  universality  property:
there is a unique such a pencil  passing  through an arbitrary second order operator acting on densities
of arbitrary weight except for three singular cases. For example, consider a second order operator $\Delta$ of weight $0$ acting on densities of a particular weight $\l$,
$\Delta\in \DD_\l(M)$,
$\Delta=S^{ab}\p_a\p_b+\dots$, where $S^{ab}$ is
a symmetric contravariant tensor field.
Then,  except for the singular cases $\l=0,{1\over 2},1$, there exists a unique canonical pencil of operators
passing through the operator $\Delta$.
These exceptional weights have deep geometric  and physical meaning. These are the same singular values that appear in the Duval--Ovsienko construction. The singularity of the  maps $\varphi_{\l\mu}$ given by equation~\eqref{isomorphism}
at $\l=0,{1\over 2},1$ is related with the existence  of non-equivalent ${\rm Diff\,}(M)$-modules (see~\cite{DuvOvs} for details).
Note that the space $\DD_{1/2}(M)$ of operators on half-densities is drastically different
from all other spaces $\DD_\l(M)$, since for the operators on half-densities there is a natural notion of a self-adjoint operator. This fact is of great importance
for the Batalin--Vilkovisky geometry (see \cite{Kh3, KhVor1}).

In the present paper, we   apply  the approach of our work~\cite{KhVor2} and  study
canonical pencils of second order operators of an arbitrary weight $\d$,
analyzing in detail the exceptional case when these operators
act on densities of weight $\lambda={1-\delta\over 2}$.
(An operator has weight $\d$ if it maps densities of weight $\l$ into  densities of weight
$\mu=\l+\d$.)   Such an operator pencil can be defined by a symmetric
contravariant tensor density
$\SS=\D{x}^\d S^{ab}\p_a\otimes \p_b$ (this field defines the principal symbol)
and a connection $\nabla$ on volume forms. Specializing the pencil to the exceptional value
  $\lambda={1-\delta\over 2}$,
  we arrive at an operator that depends only on a equivalence class of connections.
 We assign to every field  $\SS$ a  certain \emph{groupoid of connections} $C_\SS$.
  For the exceptional value  $\lambda={1-\delta\over 2}$, an operator
with the principal symbol $\SS$ depends only on an orbit of this groupoid in the space of connections.

This is particularly interesting in the case of odd symplectic structures.

Recall that  for a symplectic structure  (even or odd),  there is no distinguished affine connection associated with it (unlike Riemannian structures and the corresponding Levi-Civita connections).  On the other hand, if a symplectic structure is odd,  then the Poisson tensor is symmetric and it can be regarded as  the principal symbol $\SS$ of an odd second order differential operator or operator pencil of weight $\delta=0$.  It turns out that in spite of the absence of a distinguished affine connection, for an odd symplectic manifold there exists a distinguished class of connections on the bundle of volume forms,   such that  the connection coefficients   $\gamma_a$ (where $\nabla_a=\p_a+\gamma_a$) vanish  in some Darboux coordinates.  Connections in this distinguished class  belong to an orbit of the groupoid $C_\SS$ and the corresponding operator on half-densities is the canonical odd second order operator introduced in~\cite{Kh3}.  This operator seems to be the correct clarification of the Batalin--Vilkovisky ``odd Laplacian'' \cite{BatVyl1}.  (See~\cite{KhVor1} for details.)

This approach may be used also in the case of Riemannian geometry where the principal symbol $\SS$ is defined by the Riemannian metric.  However, in this case there exists a distinguished affine connection (the Levi-Civita connection). We would like to mention article \cite{BatBer2}, where   an interesting attempt to compare second order operators for  even Riemannian and odd symplectic structures was  made. (The similarity between even Riemannian and odd symplectic structures was pointed out in~\cite{KhVor1}.)

Another important case is a canonical pencil of operators of weight $\d=2$ on the  line.
By considering exceptional weights we arrive in particular to Schwarzian derivative.

The plan of the paper is as follows.

In the next section, we consider second order operators on the algebra of functions.
We come in this ``naive'' approach to preliminary relations between second order operators
and connections on volume forms.

In the third section, we consider first and second order operators on the algebra of   densities $\F(M)$  on a manifold $M$.
This algebra can be interpreted as a subalgebra of functions on
an auxiliary  manifold $\widehat M$. We define a canonical  invariant  scalar product on the algebra $\F(M)$.
Then we consider  the first order operators on the algebra $\F(M)$ and, in particular, the derivations of $\F(M)$ (which can be identified with   vector fields on $\widehat M$). The natural scalar product on $\F(M)$
allows to introduce a canonical divergence of the graded vector fields on $\widehat M$.
In particular,  we come to the interpretation of Lie derivatives of densities as divergence-free vector fields on the manifold $\widehat M$. After that we consider the second order operators on the algebra $\F(M)$.
Here we introduce our main construction: the self-adjoint second order operators on the algebra $\F(M)$, and consider the corresponding
operator pencils.  These considerations are due to the paper \cite{KhVor2}.

In the fourth section we consider operators of weight $\d$ acting on
densities of exceptional weight $\l={1-\d\over 2}$. For an arbitrary  contravariant symmetric tensor density $\SS$ of weight $\d$
we consider a groupoid $C_\SS$.
The orbits of the groupoid $C_\SS$ are the equivalence classes of connections such that the operators
with the principal symbol $\SS$ acting on densities of the exceptional weight $\l={1-\d\over 2}$ depend only on equivalence classes (not a choice of connection within an equivalence class). Such a groupoid was first considered in~\cite{KhVor1, KhVor2} for the case of Batalin--Vilkovisky geometry. We also give explicit description for corresponding Lie algebroids.

Then we consider various examples where these operators occur. We consider the example of operators of  weight $\d=0$  acting on half-densities on a Riemannian manifold  and on an odd symplectic supermanifold,  and the example of operators of weight $\d=2$ acting on densities of  weight $\lambda=-{1\over 2}$ on the line. In all these examples the operators depend on classes of connections on volume forms which vanish in special coordinates (such as Darboux coordinates for the symplectic case and  projective coordinates for the line).

Finally we consider the example of an odd canonical  invariant half-density introduced in~\cite{KhMkrt}. We show  that this density depends on a class of affine connections that vanish in Darboux  coordinates.

By differential operators throughout this text we mean only {\it linear} differential operators.

In this introduction and most parts of the text we speak about ``manifolds'', which can mean ordinary (or ``purely even'') manifolds, but can mean supermanifolds as well. For the simplicity of notation, we always write down the formulae as for ordinary manifolds, to avoid extra signs. However, everything extends at no extra cost to general supermanifolds. Of course, supermanifolds have to appear explicitly when we consider odd structures. For standard material on supermathematics see \cite{Berezin}, \cite{Leites} and \cite{Vor}.

A preliminary version of this text was published as~\cite{tv:hovik-chubar}.

\bigskip
\textbf{Acknowledgement.}  H.~M.~Kh.  is very happy to acknowledge the wonderful environment of the MPI Bonn, which
greatly   facilitated the work on this paper. In the course of this work, we had   many helpful conversations with V.~Yu.~Ovsienko. We are very grateful to him.

\section {Second order operators on functions}
 In what follows, $M$ is a smooth manifold or supermanifold.

Let $L=T^a(x){\p\over \p x^a}+R(x)$ be a first order operator on functions on a manifold $M$. Under a change of local coordinates $x^a=x^a(x')$, the coefficients in the operator $L$ transform as follows:
     $$
L=T^a(x){\p\over \p x^a}+R(x)=
T^a(x(x'))\,x^{a'}_a{\p\over \p x^{a'}}+R(x(x'))\,.
     \qquad\text{(Here $x^{a^\pr}_a={\p x^{a^\pr}\over \p x^a}$)}\,.
     $$
We see that  $T^a(x){\p\over \p x^a}$ is a vector field and $R(x)$ is a scalar field.

Now return to the  second order operator \eqref{secorderoperator} on a manifold $M$.
Under a change of local coordinates $x^a=x^a(x')$,
              \begin{equation}\label{secondorderoperatorfirsttime}
\Delta= {1\over 2}\left(S^{ab}(x)\p_a\p_b+T^a(x)\p_a+R(x)\right)=
    {1\over 2}\underbrace{x^{a\pr}_aS^{ab} x^{b\pr}_b}_{S^{a^\pr b^\pr}}\p_{a^\pr}\p_{b\pr}+\dots
              \end{equation}
Therefore, the top-order part  of the operator $\Delta$, ${1\over 2}S^{ab}\p_a\otimes\p_b$, defines a
symmetric contravariant tensor of  rank $2$ on $M$
(the principal symbol of the operator $\Delta= {1\over 2}\left(S^{ab}(x)\p_a\p_b+\dots\right)$.

If the tensor $S$ vanishes, then $\Delta$ becomes a first order operator, so $T^a\p_a$ is a vector field.
What about the geometric  meaning of the operator \eqref{secondorderoperatorfirsttime} in the case when the principal symbol $S\not=0$?
To answer this question, we introduce a scalar product $\langle\,\,,\,\,\rangle$ in the  space of functions on $M$
and consider the difference  of the two second order operators $\Delta^+-\Delta$,
 where $\Delta^+$ is the operator adjoint to $\Delta$ with respect to a  chosen scalar product.
A scalar product $\langle\,\,,\,\,\rangle$ on the space of functions can be defined by the following construction:
an arbitrary non-vanishing  volume form $\rh=\rho(x)|D(x)|$ on $M$ is chosen and then we set
               \begin{equation}\label{scalarproductforfunctions}
              \langle f,g\rangle_\rh:=\int_M f(x)g(x)\rho(x)|D(x)|\,.
               \end{equation}
If $x'$ are new local coordinates, so that $x^a=x^a\left(x^\pr\right)$, then in the new coordinates
the volume element $\rh$ has the form $\rho^\pr(x^\pr)|D(x^\pr)|=\rho(x)|D(x)|$:
                           $$
  \rh= \rho(x)\D{x}=\rho(x(x^\pr))\left|D(x)\over D(x^\pr)\right| |D(x^\pr)|=
 \rho(x(x^\pr))
  \left|\det\left({\p x^a\over \p x^{a^\pr}}\right)\right||D(x^\pr)|=
  \rho^\pr(x^\pr)|D(x^\pr)|\,,
                           $$
  i.e.,
                        \begin{equation*}\label{volforminnewcoordinates}
                    \rho^\pr(x^\pr)=
                      \rho(x(x^\pr))
     \left|\det\left({\p x^a\over \p x^{a^\pr}}\right)\right|.
                        \end{equation*}
In what follows, we suppose that the scalar product  is well-defined:   the manifold  $M$ is compact and orientable  and   an orienting atlas of local coordinates is chosen (all local coordinate  transformations  have positive
Jacobians: $\det \frac{\p x}{\p x^\pr} >0$\,).\footnote{A coordinate volume form $|D(x)|$ on a manifold $M$ is usually denoted by
$dx^1dx^2\dots dx^n$ or $|dx^1dx^2\dots dx^n|$. We prefer our notation $|D(x)|$ having in mind
the case of supermanifolds.}

Now return to the  operator  $\Delta$ and the adjoint  operator $\Delta^+$. For an operator $\Delta$ the operator $\Delta^+$ is defined by the relation
 $\langle\Delta f,g\rangle_\rh=\langle f,\Delta^+g\rangle_\rh$.  By integrating by parts, we obtain
         $$
\langle \Delta f,g\rangle_\rh=
     \int_M
     \underbrace
           {{1\over 2}
     \left({S^{ab}(x)\p_a\p_b}f+T^{a}(x)\p_af+R(x)f\right)
            }_{\Delta f}g(x)\rho(x)|D(x)|=
         $$
         $$
         \int_M f(x)
         \underbrace
            {
         \left(
         {1\over 2\rho}\p_a\left(\p_b\left(S^{ab}\rho g\right)\right)-
         {1\over 2\rho}\p_a\left(T^a\rho g\right)+{1\over 2}Rg
         \right)
            }_
            {\Delta^{^+} g}
        \rho(x)|D(x)|=\langle f,\Delta^+g\rangle_\rh.
         $$
The principal symbols of the operators $\Delta$ and $\Delta^+$ coincide. Thus the
difference $\Delta^+-\Delta$ is   a first order operator:
              \begin{equation}\label{vectorfieldinsymbol}
   \Delta^+-\Delta=
               \underbrace
                  {
\left( \p_bS^{ab}-T^a+ S^{ab}\p_b\log \rho \right)\p_a}_{\hbox
{vector field}}+\hbox{scalar terms}\,.
            \end{equation}
Hence, by introducing the scalar product via a choice of a volume form $\rh$ we come to the
fact that for an operator $\Delta={1\over 2}\left(S^{ab}\p_a\p_b+T^a\p_a+R\right)$
and for an arbitrary volume form $\rh=\rho(x)|D(x)|$, the expression
$\left(\p_bS^{ab}-T^a+ S^{ab}\p_b\log \rho \right)\p_a$
is a vector field.

\begin{claim}For an operator $\Delta={1\over 2}\left(S^{ab}\p_a\p_b+T^a\p_a+R\right)$\,,
 the expression
                         \begin{equation}\label{claim1}
                \gamma^a=\p_bS^{ab}-T^a
                       \end{equation}
 is an upper connection on the bundle of volume forms.
 \end{claim}

\begin{remark} The expression $\p_bS^{ab}-T^a$, up to a factor, is the so-called \emph{subprincipal symbol} of $\Delta$. See for example~\cite{hoer:volume3}.
\end{remark}
Before proving the claim,  we need to say something about connections and upper connections on  volume forms. We have collected the necessary information in the Appendix.

\begin{proof}[Proof of the claim]
Consider  the flat connection $\g^\rh\colon
\,\gamma^\rh_a=-\p_a\log\rho$\,,
defined by some chosen non-vanishing volume form  $\rh=\rho(x)|D(x)|$ (see Example
\ref{connectioncorrespondingtovolumeform}).   Since the expression
$\Y=\left(\p_b S^{ab}-T^a+ S^{ab}\p_b\log \rho \right)\p_a$ in \eqref{vectorfieldinsymbol} is a vector field
(the principal symbol of the first order operator $\Delta^+-\Delta$) and $S^{ab}\gamma^\rh_b=-S^{ab}\p_b\log\rho$
is an upper
connection, then  the sum $Y^a+S^{ab}\gamma^\rh_b$ is also an upper connection:
          $$
   S^{ab}\gamma^{\rh}_b+Y^a=-S^{ab}\p_b\log\rho+\left(\p_bS^{ab}-T^a+ S^{ab}\p_b\log \rho \right)\p_a=
         \p_bS^{ab}-T^a\,,
          $$
which proves the claim.
\end{proof}


Having in mind the above claim, we can rewrite the operator $\Delta$ on functions
in a  more convenient form:
                             $$
            \Delta f={1\over 2}\left(S^{ab}\p_a\p_b+T^a\p_a+R\right)f=
{1\over 2}\Bigl(\p_a (S^{ab}\p_bf)+L^a\p_af+Rf\Bigr)\,, \quad {\rm where}\,\,  L^a=T^a-\p_bS^{ab}\,.
                          $$
We come to the following proposition.
\begin{proposition}\label{zerothproposition}
For an arbitrary second order operator
on functions on a manifold $M$,
                  \begin{equation*}\label{firststatement}
            \Delta={1\over 2}\left(S^{ab}\p_a\p_b+T^a\p_a+R\right)=
            {1\over 2}\left(\p_a \left(S^{ab}\p_b\dots\right)+L^a\p_a+R\right)\,,
                 \end{equation*}
 the principal symbol ${1\over 2}S^{ab}$ is symmetric contravariant tensor field of rank $2$,
the subprincipal symbol $\gamma^a=-L^a=\p_bS^{ba}-T^a$ defines an upper connection on volume forms
and the function $R=2\Delta 1$ is a scalar:
              $$
    \Delta f=\frac{1}{2}\left(\p_a\left(\underbrace{S^{ab} }_{\text{\emph{tensor}}}\p_b f\right)-
    \underbrace{\gamma^a}_{\text{\emph{upper  connection}}}\p_a f+
    \underbrace{R}_{\text{\emph{scalar}}}f\right).
              $$
A second order operators on functions is fully characterized by a symmetric contravariant
tensor  of rank $2$ (the principal symbol), an upper connection on volume forms (the subprincipal symbol) and a scalar field (the value on $1$).
If the principal symbol is non-degenerate, $\det (S^{ab})\not=0$,  there arises a usual connection on volume forms:
$\gamma_a=S^{-1}_{ab}\gamma^b$.

\end{proposition}

\section{The algebra of densities and second order operators on this algebra}

\subsection{The algebra of densities $\F(M)$ and the canonical scalar product on it}

 We consider now the spaces of densities.

 As usual we suppose that $M$ is a compact orientable manifold with a chosen oriented atlas.

We  say that $\ss=s(x)\D{x}^\lambda$ is a \emph{density of weight} $\lambda$  if under a
change of local coordinates  it is multiplied by the $\l$th power of the Jacobian of the coordinate transformation:
                    $$
               \bs
           =s(x)\D{x}^\l=
 s\left(x\left(x^\pr\right)\right)\left|{Dx\over Dx'}\right|^\l\D{x^\pr}^\l=
 s\left(x\left(x^\pr\right)\right)\left|\det\der{x}{x'}\right|^\l\D{x^\pr}^\l\,.
                    $$
(A density of weight $\l=0$ is a usual function, a density of weight $\l=1$ is a  volume form.)

Denote by $\F_\lambda=\F_\l(M)$ the space of all densities of  weight $\l$ on the manifold $M$.

Denote by $\F=\F(M)=\oplus \F_\l(M)$  the space of all densities on the manifold $M$.

The space $\F_\lambda$ of densities of   weight $\lambda$ is  a vector space. It is a module over the algebra of functions on $M$.
The space $\cal F$ of all densities is itself an algebra:  If $\ss_1=s_1(x)\D{x}^{\l_1}\in \F_{\l_1}$ and
$\ss_2=s_2(x)\D{x}^{\l_2}\in \F_{\l_2}$, then their product is the
density
\begin{equation*}
    \ss_1\cdot \ss_2=s_1(x)_2(x)\,|D(x)|^{\l_1+\l_2}\in \F_{\l_1+\l_2}\,.
\end{equation*}

On the algebra $\F(M)$ of all densities  on $M$ one can consider a canonically defined scalar product $\langle\,\,,\,\,\rangle$. It is
given by the following formula: if $\ss_1=s_1(x)\D{x}^{\l_1}$ and $\ss_2=s_2(x)\D{x}^{\l_2}$, then
                 \begin{equation}\label{canonicalscalarproduct}
                   \langle \ss_1,\ss_2\rangle
                          =
                      \begin{cases}
                        \  \int_M s_1(x)s_2(x)\D{x}  \quad &\text{ if $\l_1+\l_2=1$}\,,
                                                               \\
                       \  0 \quad &\text{ if $\l_1+\l_2\neq 1$}\,.
                      \end{cases}
                 \end{equation}
(Compare this scalar product with the volume form dependent scalar product $\langle\,\,,\,\,\rangle_\rh$
on the algebra of functions
introduced in formula \eqref{scalarproductforfunctions}.)

The canonical scalar product \eqref{canonicalscalarproduct}
was introduced and  used intensively in our paper~\cite{KhVor2}. Recall briefly some  constructions from there.

The elements of the algebra $\F(M)$ are finite combinations of densities of different weights.

It is convenient to use a formal variable $t$ in   place of the coordinate volume element $\D{x}$\,.
An arbitrary density
               $
           \F\ni \ss=s_1(x)\D{x}^{\l_1}+\dots+s_k(x)\D{x}^{\l_k}
               $
can be written  as a function of $x$ and $t$ of a special form  in the variable $t$:
                         \begin{equation}\label{identification}
                      \ss=\ss(x,t)=s_1(x)t^{\l_1}+\dots+s_k(x)t^{\l_k}\,.
                         \end{equation}
For example, the density $s_1(x)+s_2(x)\D{x}^{1/2}+s_3(x)\D{x}$ can be re-written as the function
$s(x,t)=s_1(x)+s_2(x)\sqrt t + s_3(x)t$.
In what follows we will often will use this notation.

\begin{remark} With an abuse of language, we say that a function $f(x,t)$ is a  \emph{polynomial}  in $t$ if it
is a   finite sum of `monomials' of arbitrary real degrees,
$f(x,t)=\sum_{\l} f_\l(x)t^\l$, $\l\in \RR$.
(In particular, it is assumed that $t$ is an invertible variable, $t^{-1}$ makes sense.)
\end{remark}

Therefore there is a one-to-one correspondence between the elements of the algebra $\F(M)$ and the functions $s(x,t)$ polynomial in $t$.

What is the invariant meaning of the variable $t$?   The relation \eqref{identification} means
that  an arbitrary density on $M$ can be identified with a polynomial function on the `extended
manifold' $\widehat M$,
\begin{equation*}
    \widehat M=\det (TM)\setminus M\,,
\end{equation*}
the frame bundle of the determinant bundle of $M$. The natural local coordinates on $\widehat M$
induced by local coordinates $x^a$ on $M$ are $(x^a,t)$
where $t$ is the coordinate which is in  place of the volume element $\D{x}$. Note that $t\neq 0$.
Let $x^a$ and $x^{a^\pr}$ be two local coordinate systems on $M$. If $(x^a,t)$ and $(x^{a^\pr},t^\pr)$
are the corresponding induced local coordinate systems on $\widehat M$, then
                    \begin{equation}\label{transformationlaw1}
        x^{a^\pr}=x^{a^\pr}(x)\,\,{\rm and}\,\,
             t'=t\, \det\der{x'}{x}\,.
                    \end{equation}
If a function is   polynomial with respect to the local variable $t$, then
it is   polynomial with respect to the local variable $t^\pr$ as well.  (As it was mentioned before, we consider only oriented atlases, i.e., all changes of coordinates have positive determinants.)

It should be emphasized that the algebra $\F(M)$ of all densities on $M$
can be identified with the proper subalgebra   in the  algebra of functions on extended manifold $\widehat M$ consisting of all functions  that are polynomial in $t$. We  do not consider arbitrary functions of $t$.

\subsection{Derivations of the algebra    $\F(M)$}


 Consider differential operators on the algebra $\cal F$.
(We repeat that we consider only linear operators.)

Let $\X$ be a derivation of the algebra $\F$.
Then for two arbitrary densities $\ss_1$ and $\ss_2$,
                $$
          \X\left(\ss_1\cdot \ss_2\right)=\left(\X \ss_1\right)\cdot \ss_2+
          \ss_1\cdot \left(\X \ss_2\right)\,
                $$
(the Leibniz rule).  The derivations of the algebra $\F(M)$
can be identified with the vector fields on the extended manifold $\widehat M$
whose coefficients are polynomial in  $t$. We can write them as
            \begin{equation}\label{vectorfields1}
        \X=X^a(x,t){\p\over \p x^a}+X^0(x,t)\widehat \l=
\sum_\delta t^\delta \left( X^a_{(\delta)}(x){\p\over \p x^a}+X^0_{(\delta)}(x)
\widehat\l\right)\,.
                 \end{equation}
In this formula we  introduced  the Euler operator
            $$
\widehat \l=t{\p \over \p t}\,,
            $$
which is a globally defined vector field on $\widehat M$ (see the transformation law \eqref{transformationlaw1}).
The Euler operator $\widehat \l$  measures the weight   of a density:
\begin{equation*}
    \widehat \l \left(s(x)t^\l\right) =\l\, s(x)t^\l\,.
\end{equation*}
There is a natural grading on the space of  vector fields. A vector field of the form
             \begin{equation}\label{firstorderoperatorsondensities}
          \X=t^\delta\left(X^a(x)\p_a+X^0(x) \widehat \l\right)
            \end{equation}
has weight $\delta$. It transforms a density of weight $\l$ to a density of weight
$\l+\delta$.

\begin{remark}\label{vectorfield}
From now on, when speaking about vector fields on the extended manifold $\widehat M$,
we shall always suppose   that their coefficients  are polynomial in $t$, as in equation \eqref{vectorfields1}.
\end{remark}

Our next step will be to consider    adjoint operators with respect to
the canonical scalar product~\eqref{canonicalscalarproduct} on the algebra $\F$:
an  operator $\hat L^+$ is the \emph{adjoint}  to an operator $L$ if for arbitrary densities $\ss_1$ and $\ss_2$,
\begin{equation*}
    \langle \hat L \ss_1,\ss_2\rangle=\langle  \ss_1,\hat L^+\ss_2\rangle\,.
\end{equation*}
One can see that
\begin{equation*}
                    (x^a)^+=x^a\,, \quad t^+=t\,, \quad \left(\p_a\right)^+= -\p_a\, \quad
                   \text{and} \quad (\,\widehat \l\,)^+=1-\widehat \l\,.
\end{equation*}
(Here $\p_a=\lder{}{x^a}$.) Let us check the last relation. We shall write $\widehat \l^+$ for $(\,\widehat \l\,)^+$.  Let $\ss_1$ be a density of   weight $\l_1$
and $\ss_2$ be a density of   weight $\l_2$. Then $\langle \widehat \l \ss_1,\ss_2\rangle=\l_1\langle  \ss_1,\ss_2\rangle$
and   $\langle  \ss_1,(1-\widehat \l)\ss_2\rangle=(1-\l_2)\langle  \ss_1,\ss_2\rangle$. In the case when $\l_1+\l_2=1$, these scalar  products are equal since $\l_1=1-\l_2$. In the case when $\l_1+\l_2\not =1$, these scalar products both vanish (and are again equal). This proves $\widehat \l^+=1-\widehat\l$.

\begin{example} Consider  a vector field on $\widehat M$ (a derivation of the algebra of densities) $\X$:
\begin{equation*}
    \X\ss=\left(X^a(x,t)\p_a+X^0(x,t) \widehat \l\right)\ss(x,t)\,.
\end{equation*}
Then for its adjoint  operator $\X^+$ we have
\begin{equation*}
    \X^+\ss=\left[X^a(x,t)\p_a+X^0(x,t) \widehat \l\right]^+\ss=
             -\p_a\bigl(X^a(x,t)\ss \bigr) +(1-\hat\l)\bigl(X^0(x,t)\ss\bigr)\,,
\end{equation*}
so
\begin{equation*}
                  \X^+= -X^a(x,t)\p_a-X^0(x,t)\hat\l \underbrace{-\p_aX^a(x,t) + (1-\widehat \l)X^0(x,t)}\,.
\end{equation*}
Note that $\X^+$ is an operator of first order, but in general not a vector field because it contains a scalar part (underbraced).
\end{example}

\begin{definition}[Canonical divergence of vector fields on $\widehat M$]
The \emph{divergence} of a vector field  $\X$ on $\widehat M$ is defined by the formula
                    \begin{equation}\label{canonicaldivergence}
                   {\rm div\,}\X=-(\X+\X^+)=\p_aX^a+ (\widehat\l-1)X^0(x,t)\,.
                  \end{equation}
In particular, for a vector field $\X$ of  weight $\delta$,  $\X=t^\delta\left(X^a\p_a+X^0 \widehat \l\right)$
(see equation\eqref{firstorderoperatorsondensities}),
                     $$
               \div \X=t^\delta\bigl(\p_a X^a+(\delta-1)X^0\bigr)\,.
                      $$
\end{definition}
The divergence of a vector field  $\X$  vanishes if and only if it  is anti-self-adjoint
(with respect to the canonical scalar product \eqref{canonicalscalarproduct}):
             $\X=-\X^+$ $\Leftrightarrow$  $\div\X=0$.

\begin{example}\label{liedderivativeexample}
The divergence-free (= anti-self-adjoint) vector fields  of weight $\delta=0$  act on densities
as Lie derivatives. Indeed, consider a vector field $\X=X^a\p_a+X^0 \widehat \l$ of   weight $\d=0$. It defines a vector field $X=X^a\p_a$ on $M$. The condition $\div\X=\p_aX^a-X^0=0$ means that $X^0=\p_aX^a$, i.e.,
$\X=X^a\p_a+\p_a X^a \widehat \l$.  Hence for every  $\lambda$,
 $\quad \X\big\vert_{\F_\l}=X^a\p_a+\l \p_a X^a $. That means that the action
 of a divergence-free vector field $\X$ of weight $\delta=0$
 on an arbitrary density is the Lie derivative of this density with respect to the vector field $X$: for $\ss\in \F_\l$,
                     \begin{equation}\label{liederivative1}
     \quad \X \ss=(X^a\p_a+\widehat\lambda \p_a X^a)\ss={\cal L}_X \ss
        =\Bigl(X^a\p_as(x)+\l \p_a X^a s(x)\Bigr)\D{x}^\l\,.
                        \end{equation}
If $\X$ is a divergence-free vector field on $\widehat M$ of arbitrary weight,
then $\div\X=0\Leftrightarrow \X=t^\d\left(X^a\p_a+\p_aX^a{\widehat\l\over 1-\d}\right)$.
We can interpret this as a  `generalized Lie derivative': if $\delta\not=1$, then
for $\ss\in \F_\l$,
                    \begin{equation}\label{liederivative2}
                    \L_X\ss:=\X\ss=\D{x}^\d\left(X^a\p_a+\p_aX^a{\widehat\l\over 1-\d}\right)\ss=
                         \left(X^a\p_a s(x)+{\l\p_aX^a \over 1-\d}\,s\right)
                                  \D{x}^{\l+\d}\,.
                     \end{equation}
Note that $X=X^a\p_a$ in this case is a vector density on $M$ (i.e., a vector field with coefficients in densities)  of weight $\d$ , so this formula serves as a definition of $\L_X$. (One can compare this with the Nijenhuis classification of derivations of the algebra of forms and the construction of the Nijenhuis bracket of vector fields with coefficients in forms.)
\end{example}

One can consider the canonical projection $p$  of the vector fields on $\widehat M$ (the derivation of the algebra $\F(M)$) onto the vector densities on $M$. It is defined by the formula $p(\X)=\X\big\vert_{\F_0=C^{\infty}(M)}$. In coordinates, $p\colon \, \X=X^a(x,t)\p_a+X^0(x,t)\widehat\lambda \mapsto  X^a(x,t)\p_a$.  We say that a vector field is \emph{vertical} if $p(\X)=0$, i.e., if $\X=X^0(x,t)\widehat\lambda$.
The divergence of a vertical vector field $\X=X^0(x,t)\widehat\lambda$ equals to
$\div\X=(\widehat\lambda-1)X^0(x,t)$.

\begin{proposition}\label{canonicalprojection}
Let $\Pi$ be a projection of vector fields on $\widehat M$ onto the vertical vector fields
such that $\div\X=\div\left(\Pi \X\right)$. We have
           \begin{equation*}\label{projection}
  \Pi\colon\quad  \X=t^\delta\left(X^a\p_a+X^0 \widehat \l\right)\mapsto \Pi\X=
  t^\delta\left({\p_aX^a\over \delta-1}+X^0 \right)\widehat \l\,.
           \end{equation*}
 Every vector field $\X$ of  weight $\delta\not=1$
can be uniquely decomposed  into the sum of a vertical vector field and a divergence-free
vector field, which is the generalized Lie derivative \eqref{liederivative2}) with respect to the vector field $p\X$:
           \begin{equation*}
         \X=\Pi \X+\left(\X-\Pi \X\right)=\Pi \X+\L_{p\X}\,.
           \end{equation*}
\end{proposition}

One can check  the statements of this Proposition by a straightforward application of the formulae obtained above.

\smallskip

{\small What is the relation  between the canonical
divergence \eqref{canonicaldivergence} of vector fields on the extended manifold $\widehat M$ and a divergence of vector
fields on a manifold $M$ (that requires an extra structure for its definition)?
Let $\nabla$ be an arbitrary connection on volume forms.
It assigns to the vector field $\X$ on $M$ a vector field $\X_\g$ on the extended
manifold $\widehat M$ by the formula
           $\X_\g=X^a\left({\p \over \p x^a}+\widehat\lambda \gamma_a\right)$, where
$\g=\{\gamma_a\}$ is the connection form for $\nabla$ in coordinates $x$,
 ($\nabla_a\D{x}=\gamma_a\D{x}$).
A connection $\nabla$ defines a divergence of vector fields on $M$
via the canonical divergence \eqref{canonicaldivergence} on  $\widehat M$: for every vector field $\X$ on the manifold $M$,
       \begin{equation}\label{divergenceonmanifold}
      {\div}_\g\X:=\div\X_\g=\left({\p X^a\over \p x^a}-\gamma_a X^a\right)\,.
       \end{equation}
A non-vanishing volume form $\rh=\rho(x)\D{x}$ defines the flat connection $\gamma^\rh_a=-\p_a\log \rho$
(see equation \eqref{exampleofflatconnection} and example \ref{connectioncorrespondingtovolumeform}).
The formula \eqref{divergenceonmanifold} gives the familiar formula
(see also equation \eqref{firstequation}) for the divergence of vector fields on
a manifold equipped with a volume form:
        \begin{equation}\label{divergenceonmanifoldwithvolumeform}
      {\div}_\rh\X:=\div\X_{\g^\rh}=\left({\p X^a\over \p x^a}+ X^a\p_a\log\rho\right)
       ={1\over \rho}{\p\over \p x^a}\left(\rho X^a\right)\,.
       \end{equation}
If we consider the connection on volume forms corresponding to an affine connection on $M$ (see Example \ref{secondexample}),
we come to ${\rm div\,}_{ \nabla}\X=\nabla_aX^a=
(\p_a X^a+X^a\Gamma_{ab}^{b})$.  On a Riemannian manifold $M$, the
Riemannian metric defines the connection on volume forms $\gamma_a=-\p_a\log \sqrt {\det g}$
(via the Levi-Civita connection or via the invariant volume element  $\rh_g$).    We come to the familiar formula
          \begin{equation*}\label{divergenceonRiemannian}
         {\div}_g\X=\left({\p X^a\over \p x^a}+ X^a\p_a\log\sqrt{\det  g}\right)=
          {1\over \sqrt {\det g}}{\p\over \p x^a}\left(\sqrt {\det g} X^a\right)\,.
          \end{equation*}
defining the Riemannian divergence of vector fields.

}

\subsection{Second order operators   on the algebra $\F(M)$}

Let us now turn to differential operators of order $\leq 2$ on the algebra of densities.

First of all, a general remark about the definition of the $n$th order operators.
A $0$th order operator on the algebra $\F(M)$ is just a multiplication operator (the multiplication by a given density).
A linear operator $L$   on the algebra $\F(M)$ is a differential operator of order $\leq n$ (or an $n$th order operator) if for an arbitrary $\ss\in \F(M)$
the commutator with the multiplication operator $[L,\ss]=L\circ \ss-\ss\circ L$ is an operator of order $\leq n-1$.

One can see that if $L$ is a differential operator on  $\F(M)$ of order $\leq n$, then the operator
$L+(-1)^nL^+$ is also of order $\leq  n$ and the operator $L-(-1)^nL^+$ is
of order $\leq n-1$.   We have the following statement.

\begin{proposition}\label{decomposition}
An arbitrary $\mathrm n${th} order operator can be canonically decomposed into  the sum
of a self-adjoint and an anti-self-adjoint operators:
     \begin{equation*}\label{decompositionof n-th order operator}
     L=\underbrace{{1\over 2}\left(L+(-1)^nL^+\right)}_{\text{\emph{operator of order $n$}}}
       +\underbrace{{1\over 2}\left(L-(-1)^nL^+\right)}_{\text{\emph{operator of order $\leq n-1$}}}\,.
 \end{equation*}
An operator of even order $n=2k$ is the sum of a self-adjoint operator of order $2k$
and  a anti-self-adjoint operator of  order $\leq 2k-1$. An
operator of odd order $n=2k+1$ is the sum of a anti-self-adjoint operator of  order $2k+1$
and  a self-adjoint operator of order $\leq 2k$.
\end{proposition}

Operators of order $0$ are evidently self-adjoint.

Let $L=\X+B$ be a first order anti-self-adjoint operator, where $\X$ is a vector field on $\widehat M$
and $B$ is a scalar term
(density). We have $L+L^+=0=\X+\X^++2B=0$. Hence $L=\X+{1\over 2}\div\X$.


Now let us study the self-adjoint second order operators on the algebra of densities $\F(M)$.
Let $\Delta$ be a second order operator of  weight $\delta$ on $\F(M)$.
In local coordinates,
            \begin{equation}\label{secondordergeneral}
   \Delta={t^\delta\over 2}
            \left(
    \underbrace{S^{ab}(x)\p_a\p_b+\widehat\lambda B^a(x)\p_a+\widehat \l^2 C(x)}_{\hbox{second order derivatives}}+
\underbrace{D^a(x)\p_a+\widehat\lambda E(x)}_{\hbox{first order derivatives}}+
          F(x)
             \right)\,.
           \end{equation}
Impose a normalization condition
\begin{equation}\label{normalisationcondition}
\Delta (1)=0\,,
\end{equation}
so that the density $F\D{x}^\delta$ in \eqref{secondordergeneral} vanishes. The  operator  $\Delta^+$ adjoint to $\Delta$ equals
\begin{multline}
                \Delta^+=\frac{1}{2}
            \Biggl(
            \p_b\p_a
            \left(S^{ab}t^\d\dots\right)-\p_a\left(B^a \widehat \l^+ (t^\d\dots)\right)+
            (\widehat \l^+)^2 \left(C t^\d\dots\right)-
            \p_a\left(D^a t^\d\dots\right)+E \widehat \l^+ \left(t^\d\dots\right)
            \Biggl)=
                \\
                {t^\delta\over 2}
            \left(
            S^{ab}\p_a\p_b+2\p_bS^{ba}\p_a+\p_a\p_bS^{ba}\right)+
                     \\
                     {t^\delta\over 2}
            \Biggl(
            \left(\widehat\lambda+\delta-1\right)\left(B^a\p_a+\p_bB^b\right)+
            \left(\widehat\lambda+\delta-1\right)^2C-
            \left(\widehat\lambda+\delta-1\right)E-D^a\p_a-\p_bD^b
            \Biggr)\,.
\end{multline}
Comparing this operator with operator \eqref{secondordergeneral}
    we see that the condition $\Delta^+=\Delta$ implies that
             \begin{equation}\label{selfadjointoperator0}
   \Delta={t^\delta\over 2}
            \Biggl(
    S^{ab}\p_a\p_b+\p_bS^{ba}\p_a+\left(2\widehat\lambda+\delta-1\right) \gamma^a\p_a+\widehat\lambda\p_a\gamma^a+
         \widehat\lambda \left(\widehat\lambda+\delta-1\right)\theta
             \Biggr)\,.
           \end{equation}
Here for convenience we have denoted   $\gamma^a(x)=2B^a(x)$ and $\theta(x)=C(x)$.
By studying how the coefficients of the operator transform under a change of coordinates, we come to the following statement.
 \begin{proposition}[See \cite{KhVor2}]\label{theorem1}
Let $\Delta$ be an arbitrary linear second order self-adjoint operator of weight $\delta$
on the algebra of densities $\F(M)$ normalized by the condition $\Delta(1)=0$.  Then in arbitrary local coordinates this operator
has the form \eqref{selfadjointoperator0}. The coefficients of this operator have the following
geometric meaning:
\begin{itemize}
    \item   $S^{ab}(x)$ are components of a symmetric
    contravariant tensor density of  weight $\delta$. Under a change of local coordinates
         $x^{a^\pr}=x^{a^\pr}(x)$ they transform in the following way:
                           $$
                    S^{a^\pr b^\pr}=J^{-\delta}
                  x^{a^\pr}_ax^{b^\pr}_b S^{ab}\,,
                   $$
        \item   $\gamma^a$ are coefficients  of an upper connection-density of
         weight $\delta$ (see \eqref{transformationofupperconnection} above).
Under a change of local coordinates
 $x^{a^\pr}=x^{a^\pr}(x)$  they transform in the following way:
                           $$
                    \gamma^{a^\pr}=J^{-\delta}
                  x^{a^\pr}_a\left(\gamma^a+S^{ab}\p_b \log J\right)\,,
                   $$
\item $\theta$ transforms in the following way:
                               $$
                    \theta'=J^{-\delta}
        \left(\theta+2\gamma^a\p_a \log J+\p_a\log J\, S^{ab}\p_b \log J\right)\,.
                               $$
\end{itemize}

Here $J=\det\left({\p x'\over \p x}\right)$  and $x^{a^\pr}_a$ are shorthand notations
for the derivatives: $x^{a^\pr}_a=\lder{x^{a^\pr}\!}{x^a}$.
\end{proposition}

We call the object  $\theta(x)\D{x}^\d$ a \emph{Brans-Dicke function}
\footnote{Its transformation is similar to the transformation
of the Brans-Dicke ``scalar" $g^{55}$  in the Kaluza-Klein reduction of the $5$-dimensional gravity to gravity+electromagnetism.}.

\begin{corollary}[From Proposition~\ref{theorem1}]\label{remarkaftertheorem2aboutnondegeneratesymbol}
A given tensor density
${\bf S}=S^{ab}\D{x}^\delta\,\p_a\otimes\p_b$ of weight $\d$ and a connection
on volume forms $\g$  uniquely define a second order self-adjoint operator \eqref{selfadjointoperator0} with $\mathbf S$ as the principal symbol, and with the upper connection $\gamma^a=S^{ab}\gamma_b$
and the Brans-Dicke function $\theta(x)=\gamma_aS^{ab}\gamma_b$.
We denote this operator $\Delta(\SS,\g)$.

The converse implication holds
if the principal symbol $\SS$ is non-degenerate: a second order self-adjoint operator $\Delta$
of weight $\delta$  with a non-degenerate principal symbol $\SS$, obeying the normalization condition \eqref{normalisationcondition}, uniquely defines a connection on volume forms $\g$ such that
$\Delta=\Delta (\SS,\g)+\hat\l(\hat\l+\delta-1)F$, where $F$ is a density of weight $\delta$, so that the
Brans-Dicke function $\theta$
is  $\theta=\gamma_a\gamma^a+F=\gamma_aS^{ab}\gamma_b+F$.
\end{corollary}

\begin{remark}\label{remarkaftertheorem1}
Let $\Delta$ be the self-adjoint operator \eqref{selfadjointoperator0} and $\g^{\pr}=\{\gamma_a^{\pr}\}$ be
 an arbitrary connection on volume forms: $\nabla\D{x}=\gamma_a^{\pr}\D{x}$. Then
for an upper connection-density  in the equation \eqref{selfadjointoperator0} the difference
$(\gamma^a-S^{ab}\gamma_b^{\pr}) \D{x}^\d$ is a vector density of  weight $\delta$.

The difference between the operator $\Delta$ and the operator
$\Delta(\SS,\g^\pr)$ is an operator which can be expressed via
a generalized Lie derivative~\eqref{liederivative2}
and a density of weight $\delta$.
\end{remark}

Let us consider examples.

First consider an example of an operator \eqref{selfadjointoperator0} with a
degenerate principal symbol $S^{ab}\D{x}^\delta$.

\begin{example}\label{operatorwothdegeneratesymbol}
Let $X=X^a{\p \over \p x^a}$ and $Y=Y^a{\p \over \p x^a}$ be two vector fields on the manifold $M$.
Recall the operator  of the Lie derivative $\L_X=X^a\p_a+\widehat\lambda \p_aX^a$ (see equation \eqref{liederivative1}). Consider the operator
             \begin{equation*}\label{operatordegenerate}
                \Delta={1\over 2}\left(\L_X\L_Y+\L_Y\L_X\right)=
      {1\over 2}\left(X^a\p_a+\widehat\lambda \p_aX^a\right)\left(Y^b\p_b+\widehat\lambda \p_bY^b\right)+
      \left(X\leftrightarrow Y\right)\,.
                  \end{equation*}
It is a self-adjoint operator because a Lie derivative is a anti-self-adjoint operator. Calculating this operator
and comparing  it with the expression
\eqref{selfadjointoperator0} we come to
                $$
            S^{ab}=X^aY^b+Y^bX^a\,, \quad \gamma^a=\left(\p_bX^b\right)Y^a+\left(\p_bY^b\right)X^a\,, \quad
\theta=\left(\p_aX^a\right)\left(\p_bY^b\right)\,.
                $$
We see that in  general   (if the dimension of the manifold is greater than $2$) this operator has a degenerate principal symbol
and the upper connection $\gamma^a$ does not uniquely define  a genuine connection.
\end{example}

\subsection{Canonical pencils}
Note that a differential operator $L$ on the algebra  of densities $\F(M)$ defines a pencil $\{L_\l\}$ of operators on spaces $\F_\l$:
    $L_\l=L\big\vert_{\F_\l}$.
A self-adjoint operator of second order $\Delta$ on the algebra of densities (see equation
        \eqref{selfadjointoperator0})
        defines an operator pencil $\{\Delta_\l\}, \l\in \RR$, where
                              $$
              \Delta_\l=\Delta\big\vert_{\F_\l}=
                              $$
             \begin{equation}\label{canonicalpencil}
                             =
             {t^\delta\over 2}
            \Bigl(
    S^{ab}(x)\p_a\p_b+\p_bS^{ba}\p_a+\left(2\l+\delta-1\right) \gamma^a(x)\p_a+\l\p_a\gamma^a(x)+
         \l \left(\l+\delta-1\right)\theta(x)
             \Bigr)\,.
           \end{equation}
This pencil is defined by a symmetric tensor density  $\SS^{ab}=S^{ab}(x)\D{x}^\d$, an
upper connection   $\gamma^a$ and a Brans--Dicke function  $\theta(x)$. Respectively, a self-adjoint operator $\Delta(\SS,\g)$ on
the algebra of densities defined by a tensor density  $\SS=S^{ab}(x)\D{x}^\d$
and a genuine connection $\g$  (see Corollary \ref{remarkaftertheorem2aboutnondegeneratesymbol})
defines an operator pencil that we denote $\Delta_\l(\SS,\g)$, with the Brans--Dicke function  $\theta(x)=\gamma_aS^{ab}\gamma_b$.

An arbitrary operator $\Delta_\l$ of   weight $\delta$ maps densities of weight $\l$ to densities
of weight $\l+\d$.
Its adjoint   $(\Delta_\l)^+$
maps densities of weight $1-\l-\delta$ to densities
of weight $1-\l$. The condition $\Delta=\Delta^+$ of the self-adjointness for an operator $\Delta$ on the algebra    $\F(M)$ is equivalent to the condition
                                \begin{equation}\label{selfadjointforpencil}
                                (\Delta_\l)^+=\Delta_{1-\l-\d}\,
                                \end{equation}
for the corresponding pencil. We shall refer to operator pencils satisfying this condition  as to the \emph{self-adjoint pencils}. The condition $\Delta(1)=0$ for operators on $\F(M)$ becomes $\Delta_0(1)=0$ for the pencils. We shall refer to the pencils satisfying $\Delta_0(1)=0$, as to \emph{normalized}.

\begin{example}\label{beltrami-laplaceoperatorondensities2}    Let $\rh=\rho(x)\D{x}$ be an non-vanishing volume form on a Riemannian manifold $M$.
Consider an operator
on functions $\Delta$  defined by $\Delta f= \frac{1}{2}\,\div\grad f=
\frac{1}{2}\,\frac{1}{\rho}\,\p_a (\rho g^{ab}\p_b{f})$, see   equations \eqref{firstequation} and \eqref{divergenceonmanifoldwithvolumeform}.
(If $\rh=\sqrt {\det g}\D{x}$, this is just the Laplace-Beltrami operator \eqref{laplacebeltrami}.)
Using the operator $\Delta$ on functions, we may introduce an operator pencil
                 \begin{equation*}
   \Delta_\l=\rh^\l\circ \Delta\circ {1\over \rh^\l}\,,
                \end{equation*}
so that for $\ss\in \F_\l$,
\begin{equation*}
      \Delta_\l\ss=  \rh^\l\div\grad \bigl(\rh^{-\l}\ss\bigr)\,,
\end{equation*}
where $\div={\div}_{\rh}$.
One can see that this pencil corresponds to a self-adjoint operator
(see relation \eqref{selfadjointforpencil}). It coincides with the canonical
pencil \eqref{canonicalpencil} of weight $\delta=0$
in which the principal symbol is given by the Riemannian metric,
$S^{ab}=g^{ab}$,   the connection is a flat connection defined by the volume element $\rh$
(see formula \eqref{exampleofflatconnection} and Example \ref{connectioncorrespondingtovolumeform}),
$\gamma^a=-g^{ab}\p_b\log \rho$, and $\theta=\gamma^a\gamma_a$.
\end{example}

The canonical pencil \eqref{canonicalpencil} has many interesting properties (see \cite{KhVor2} for details).
In particular, it has the following ``universality'' property, directly following from equation~\eqref{canonicalpencil}:

\begin{theorem}\label{coroftheorem}
Let
$\Delta\colon \F_{\l_0}\to \F_{\mu_0}$  be an arbitrary  second order differential operator of weight $\delta$
acting on the space $\F_{\l_0}$ of densities of weight $\l_0$, $\mu_0=\l_0+\d$.
If   $\l_0\not=0, \mu_0\not=1$ and $\l_0+\mu_0\not=1$,  there exists a unique normalized self-adjoint pencil $\Delta_{\l}$ given by~\eqref{canonicalpencil}
that passes through the operator $\Delta$, i.e., $\Delta_{\l}|_{\l=\l_0}=\Delta$.
If the operator $\Delta$ is given by an expression
$\Delta=A^{ab}\p_a\p_b+A^a\p_a+A(x)$, then the relations
    \begin{equation*}\label{definingcanonicalpencil}
    \begin{cases}
   {1\over 2} S^{ab}&=A^{ab}\,,\cr
   {1\over 2}\left((2\l_0+\d-1) \gamma^a+\p_bS^{ba}\right) &=A^a\,,\cr
   {1\over 2}\left(\l_0\p_a\gamma^a+ \l_0(\l_0+\d-1)\theta\right) &=A\,,\cr
       \end{cases}\qquad (\l_0\not=0, \l_0+\mu_0\not=1, \mu_0\not=1)\,.
    \end{equation*}
uniquely define the principal symbol, the upper connection and the Brans-Dicke field that uniquely define the pencil $\Delta_{\l}$ by~\eqref{canonicalpencil}.
\end{theorem}

(To emphasize: $\l_0$ and $\delta$ in the above theorem are both concrete numbers, while in the formulas for the pencil $\Delta_{\l}$, $\l$ is a parameter, which can be specified in particular to $\l=\l_0$.)

The ``universality'' property provides a beautiful interpretation of the canonical map
$\varphi_{\l\mu}$ in the relation \eqref{isomorphism}. For convenience, let us speak about $\l_0$ and $\mu_0$.
According to the Theorem,  we can ``draw"  a unique canonical pencil  $\Delta_{\l}$ through an arbitrary
operator $\Delta_{(\l_0)}$
acting on densities of weight $\l_0$.  Then the map $\varphi_{\l_0\mu_0}$ maps this operator $\Delta_{(\l_0)}$ to
an operator $\Delta_{(\mu_0)}$ acting on densities of weight $\mu_0$ obtained by the specialization of the pencil $\Delta_{\l}$ to the value $\l=\mu_0$.

\section {Operators on a manifold depending on a class  of connections}

In this section,  we will return from differential operators  on the algebra $\F(M)$   to operators on a manifold $M$ acting on densities of  particular weight.

\subsection {Operators of weight $\d$ acting on densities of weight $1-\d\over 2$}

Theorem~\ref{coroftheorem} states that for a second order operator $\Delta\colon \F_\l\to \F_\mu$ of weight $\delta$,
for all values of weight $\l$ except for the cases $\l=0$, $\mu=1$ or $\l+\mu=1$,
there is a unique canonical pencil \eqref{canonicalpencil}
which passes through the operator $\Delta$. Consider now an exceptional case when an operator $\Delta\colon \F_\l\to \F_\mu$ is such that  $\l+\mu=1$, i.e., $\Delta$ acts on densities of  weight $\l={1-\d\over 2}$ and maps them to densities of weight $\mu={1+\d\over 2}$.

Suppose a canonical pencil~\eqref{canonicalpencil} is   still given  and  consider  its specialization to the singular value $\l=\frac{1}{2}(1+\delta)$. Denote this specialization by $\Delta_{\text{sing}}$, so $\Delta_{\text{sing}}$ is an operator $\F_{{1-\delta\over 2}}\to \F_{{1+\delta\over 2}}$. We have
             $$
             \Delta_{\text{sing}}=
\left(\Delta_\l\right)\big\vert_{\l={1-\d\over 2}}=
             {t^\delta\over 2}
            \left(
    S^{ab}(x)\p_a\p_b+\p_bS^{ba}\p_a+\l\p_a\gamma^a(x)+
         \l \left(\l+\delta-1\right)\theta(x)
             \right)=
            $$
\begin{equation}\label{singvalueofcanpencil}
             ={\D{x}^\delta\over 2}
            \Biggl(
    S^{ab}\p_a\p_b+\p_bS^{ba}\p_a+
              {1-\d\over 2}
              \left(
              \p_a\gamma^a+
                  \frac{\d-1}{2}\,\theta
             \right)
             \Biggr)\,.
\end{equation}

On the other hand, consider     an arbitrary second order differential operator $\Delta$
of weight $\d$  which  acts on densities of weight $1-\d\over 2$, so that
  $\Delta\colon \F_{{1-\delta\over 2}}\to \F_{{1+\delta\over 2}}$.
We shall compare this operator with the operator $\Delta_\sing$.
  The operator
$\Delta^+$ which is adjoint to the operator $\Delta$ also acts from the
space $\F_{{1-\delta\over 2}}$ into the space $\F_{{1+\delta\over 2}}$,
since $\l+\mu={1-\delta\over 2}+{1+\delta\over 2}=1$ (compare with formula \eqref{selfadjointforpencil}).
Hence the operator $\Delta$ can be uniquely decomposed into the sum of a second order
self-adjoint operator and an anti-self-adjoint operator of   order $\leq 1$.
This anti-self-adjoint operator is just the generalized Lie derivative \eqref{liederivative2} along a vector density $X$ of weight $\delta$ specialized to $\l={1-\d\over 2}$:
                             $$
\Delta^+-\Delta=\L_X\big\vert_{\F_{1-\d\over 2}}=
\D{x}^\d\left(X^a\p_a+{1\over 2}\,\p_a X^a\right)\,.
                   $$
Operator $\Delta_{\rm sing\,}$ in formula \eqref{singvalueofcanpencil}
belongs to canonical pencil, it
is a self-adjoint operator: $\Delta_\sing^+=\Delta_\sing$.
The difference of two self-adjoint operators of second order
with the same principal symbol is a self-adjoint operator of order $\leq 1$.
Hence it is a zeroth order operator of multiplication by a  density.
These considerations imply the following statement:

\begin{corollary}[From Proposition~\ref{theorem1}]\label{corollary2}
Let $\Delta$ be an arbitrary second order operator of weight $\delta$
acting on the space of densities of weight ${1-\d\over 2}$, with the principal symbol $\SS^{ab}=S^{ab}\D{x}^\d$.
Let $\Delta_{\text{\emph{sing}}}$ be  an operator of the same weight $\delta$ acting on densities of weight ${1-\d\over 2}$ and belonging to an {arbitrary}
canonical pencil \eqref{canonicalpencil},  with the same principal symbol $\SS^{ab}=S^{ab}\D{x}^\d$.

Then the difference $\Delta-\Delta_{\text{\emph{sing}}}$ is an operator of order $\leq 1$ which is the sum of a generalized Lie derivative \eqref{liederivative2} with respect to a vector density $X$ and a zeroth order operator of multiplication by a density $\FF=F(x)\D{x}^\delta$:
               $$
       \Delta=\Delta_{\text{\emph{sing}}}+\L_X+\FF\,.
               $$
If the operator $\Delta$ is self-adjoint, $\Delta^+=\Delta$, then $X=0$ and the generalized Lie derivative vanishes from the formula.
\end{corollary}

Let us fix a pencil~ \eqref{canonicalpencil}.
It follows from this Corollary  that if
the operator $\Delta\colon \F_{1-\d\over 2}\to \F_{1+\d\over 2}$ with the principal symbol $\SS$ is a self-adjoint operator, then it is given in local coordinates by the expression
         \begin{equation*}
             \Delta=
             {1\over 2}\left(
    S^{ab}(x)\p_a\p_b+\p_bS^{ba}(x)\p_a+
               U_\SS(x)
             \right)\D{x}^\delta \,,
             \end{equation*}
where
               \begin{equation*}\label{appearanceofsingoperator}
            U_\SS(x)\D{x}^\delta=
              {1-\d\over 2}\,
              \left(
              \p_a\gamma^a(x)+
                  {\d-1\over 2}\,\theta(x)
             \right)
              \D{x}^\delta+F(x)\D{x}^\d\,.
               \end{equation*}
Here  $\gamma^a,\theta$ are the upper connection and the Brans-Dicke field defining the
pencil \eqref{canonicalpencil}, and $F(x)\D{x}^\d$ is some density.
In particular,  the self-adjoint operator  $\Delta\colon \F_{1-\d\over 2}\to \F_{1+\d\over 2}$
belongs to the canonical pencil defined by the same principal symbol $\SS$  and the upper connection $\gamma^a$
but a possibly different $\theta^\pr=\theta-{4F\over (\d-1)^2}$.  It may belong to many other pencils with
different upper connections.

A self-adjoint operator $\Delta$ acting on densities of the exceptional weight $\l={1-\d\over 2}$ \emph{does not} define uniquely  a canonical pencil  to which it belongs.  Therefore we arrive to the idea of a certain groupoid acting on connections (for a given `principal symbol' $\SS$).

\subsection {The groupoid of connections $C_\SS$}

We shall define now a certain `groupoid of connections' $C_\SS$ associated with a contravariant tensor density  $\SS^{ab}=S^{ab}\D{x}^\d$ of weight $\delta$.

Consider the space $\A$ of all connections (covariant derivatives) on volume forms  on a manifold $M$.
 This is an affine space associated with the vector space of covector fields on $M$:
the difference of two connections $\nabla$ and $\nabla'$ is a covector field (or differential $1$-form):
\begin{equation*}
\nabla-\nabla'=\g-\g'=\X=X_adx^a \quad \text{where} \quad X_a=\gamma_a-\gamma_a'\,.
\end{equation*}
(Henceforth we shall be using boldface letters for denoting vector or covector fields on $M$.)


We first consider a trivial groupoid whose set of points is the space of connections $\A$ and the set of arrows is the set $ \left\{\g{\buildrel \X\over \longrightarrow} \,\g^\pr\right\}$, where $\g,\g^\pr\in \A$ are connections and $\g^\pr=\g+\X$,
so that  $\X$ is the difference of connections, a covector field.
We have
      \begin{equation}\label{trivialabeliangroupid}
      -\left(\g_1{\buildrel \X\over \longrightarrow} \,\g_2\right)= \g^\pr{\buildrel -\X\over \longrightarrow}\, \g\,,\qquad
 \g_1{\buildrel \X\over \longrightarrow}\, \g_2\ + \
               \g_2{\buildrel \Y\over \longrightarrow} \,\g_3 =
      \gamma_1{\buildrel \X+\Y\over \longrightarrow}\, \gamma_3\,.
       \end{equation}

The groupoid of connections $C_\SS$ that we want is introduced as a subgroupoid of this trivial groupoid. This is done as follows.

Pick an arbitrary contravariant symmetric tensor density of weight $\delta$:
  $\SS^{ab}(x)=S^{ab}(x)\D{x}^\d$. The tensor density $\SS$ and an arbitrary connection $\g$ on volume forms
define the  self-adjoint operator $\Delta(\SS,\g)$  on the algebra of densities.
It is the operator defined in the equation \eqref{selfadjointoperator0}. Here the principal symbol is $\SS$, the upper connection
is $\gamma^a=S^{ab}\gamma_b$ and the Brans-Dicke function is $\theta=\gamma_a\gamma^a$
(see  Corollary \ref{remarkaftertheorem2aboutnondegeneratesymbol}).
Consider the corresponding operator pencil and the operator $\Delta_\sing(\SS,\g)$ which belongs to this pencil
and  acts on densities of weight $1-\d\over 2$:
\begin{multline}
    \Delta_{\rm sing\,}(\SS,\g)=
            \Delta\left({\SS,\g}\right)\big\vert_{\F_{1-\delta\over 2}}=\\
            \label{operatorsingular2}
             {\D{x}^\delta\over 2}
            \left(
    S^{ab}\p_a\p_b+\p_bS^{ba}\p_a+
          {1-\d\over 2}
          \left(
          \p_a\gamma^a+
              {\d-1\over 2}\,\gamma_a\gamma^a
             \right)
             \right)\,.
\end{multline}
Thus an arbitrary contravariant symmetric tensor density
 $\SS^{ab}$  of weight $\d$ and an arbitrary connection  on forms $\g$ define a  self-adjoint operator
$\Delta_{\rm sing\,}(\SS,\g)$ by equation~\eqref{operatorsingular2}.
The ``pseudoscalar'' part of this operator is equal to
             \begin{equation}\label{pseudoscalarpart}
            U_{\SS,\g}(x)\D{x}^\d =
               \left(
             {1-\d\over 2}
             \right)
                \left(
          \p_a\gamma^a+
              {\d-1\over 2}\gamma_a\gamma^a
             \right)
             {\D{x}^\delta\over 2}\,.
                \end{equation}
Let $\g$ and $\g'$ be two different connections. The difference of the two operators
$\Delta_\sing(\SS,\g)$ and $\Delta_\sing(\SS,\g')$ with the same principal symbol
$\SS^{ab}=S^{ab}(x)\D{x}^\d$ is a scalar density of weight $\delta$.   Let us calculate this density.
If $\g'=\g+\X$,  then
\begin{multline}
\Delta_\sing(\SS,\g')-\Delta_\sing(\SS,\g)=U_{\SS,\g'}(x)\D{x}^\d-U_{\SS,\g}(x)\D{x}^\d=\\
                                              \left(
                       {1-\d\over 4}
                         \right)
                         \left(
                             \p_a\gamma^{\pr a}+
              {\d-1\over 2}\gamma'_a\gamma^{\pr a}-
              \p_a\gamma^a-
              {\d-1\over 2}\gamma_a\gamma^a\right)\D{x}^\d=\\
                                 \left(
                              {1-\d\over 4}
                                 \right)
 \left(\p_a(S^{ab}X_b)+(\d-1)\gamma_a (S^{ab}X_b) +{\d-1\over 2}X_aS^{ab}X_b\right)\D{x}^\delta=\\
                                 \label{differenceoftwosingularoperators}
                                  {1-\d\over 4}
                                  \left(
                    {\div}_{\g} \X+{\d-1\over 2}\,\X^2
                                   \right)\,.
 \end{multline}
 Here $\div_{\g}\X$ is the divergence of a vector density $\X$ on $M$ with respect to a connection $\g$
 (see \eqref{divergenceonmanifold}) and the scalar square $\X^2$ is defined with the help of the tensor density $\SS$.
 With some abuse of notation we denote a covector field $X_adx^a$ and the vector density $X^a\D{x}^\d=S^{ab}X_b\D{x}^\delta$ of weight
 $\d$  by the same letter $\X$.


Now we can define our groupoid.

\begin{definition}
Let $\SS^{ab}=S^{ab}(x)\D{x}^\d$ be a contravariant symmetric tensor density of weight $\delta$.
The \emph{groupoid of connections} $C_\SS$ associated with $\SS^{ab}$ is a subgroupoid of the trivial groupoid defined above
with the same set of objects (connections on volume forms) and the subset of arrows $\left\{{\g{\buildrel \X\over \longrightarrow}
  \g^\pr}\right\}$ specified by the condition  that the operators $\Delta_\sing(\SS,\g)$  and $\Delta_\sing (\SS,\g')$
 defined by  formula \eqref{operatorsingular2} coincide:
          \begin{equation}\label{definitionofgroupoid}
 C_\SS=\left\{{\g{\buildrel \X\over \longrightarrow}
 \,\g^\pr} \ \Bigl|\Bigr. \  \Delta_\sing(\SS,\g^\pr)=\Delta_\sing(\SS,\g)\right\}.
 \end{equation}
\end{definition}
In  other words, an arrow $\grarrow$ belongs to the groupoid $C_\SS$ if the two canonical pencils $\Delta_\l(\SS,\g)$
and $\Delta_\l(\SS,\g^\pr)$ intersect at the operator $\Delta_\sing (\SS,\g)$.
Using formula \eqref{differenceoftwosingularoperators} for the difference of operators  $\Delta_\sing(\SS,\g')$ and  $\Delta_\sing(\SS,\g)$, we can
rewrite the definition \eqref{definitionofgroupoid} of the groupoid $C_\SS$ in the following equivalent way:
\begin{equation}\label{conditionofgroupid2}
 C_\SS=
 \left\{{\g{\buildrel \X\over \longrightarrow}
 \,\g^\pr} \ \Bigl|\Bigr. \  {\div}_{\g} \X+\frac{\d-1}{2}\,\X^2=0 \right\}\,.
\end{equation}

(We consider the case $\d\not=1$. The case $\delta=1$ is trivial.\footnote{ In this case   the operators
$\Delta_\sing(\SS,\g)$ do not depend on a connection $\g$ at all. The principal symbol
$\SS^{ab}=S^{ab}\D{x}$ defines a canonical  operator
$\Delta(\SS)\colon \F_0\to \F_1$ such that in local coordinates
$\Delta(\SS) f=\p_a\left(S^{ab}\p_b f\right)\D{x}$.   The groupoid $C_\SS$ is here the trivial
groupoid of all connections.})


The groupoid $C_\SS$ defines a partition of the space of connections $\A$ into orbits (or equivalence classes).
Denote by  $[\g]$ the orbit of  a connection $\g$\,:
\begin{equation*}\label{orbitofelementinabeliangroupoid}
          [\g]=\{\g'\colon \quad\grarrow\in C_\SS\}\,.
\end{equation*}

We may summarize our constructions in the following proposition.
\begin{proposition}
An arbitrary contravariant symmetric tensor density $\SS^{ab}=S^{ab}(x)\D{x}^\d$ of  weight $\d$
defines the groupoid of connections $C_\SS$ and
a family of second order differential operators
of weight $\delta$   acting on densities of   weight $1-\d\over 2$:
                         $$
               \Delta([\g])=\Delta_{\text{\rm sing}}(\SS,\g)\colon \F_{1-\d\over 2}\to \F_{1+\d\over 2}\,.
                $$
The operators in this family have the same principal symbol $\SS$ and they depend on equivalence classes of connections that  are the orbits
of the groupoid $C_\SS$.
\end{proposition}

\begin{remark}\label{cocycleconditionforbv0}
      Let $\g_1,\g_2$ and $\g_3$ be three arbitrary connections.
Consider the corresponding arrows
$\g_1{\buildrel \X\over \longrightarrow} \g_2$, $\g_2{\buildrel \Y\over \longrightarrow} \g_3$
and $\g_1{\buildrel \X+\Y\over \longrightarrow} \g_3$.
 We have
 $\g_1{\buildrel \X\over \longrightarrow}\g_2+\g_2{\buildrel \Y\over \longrightarrow} \g_3=
        \g_1{\buildrel \X+\Y\over \longrightarrow} \g_3$.
This means that for the
non-linear differential equation  ${\div}_{\g} \X+{\d-1\over 2}\,\X^2=0$ the following property holds:
                \begin{equation*}\label{symmetryofbvequation}
       \begin{cases}
       {\div}_{\g_1} \X+{\d-1\over 2}\,\X^2=0\cr
       {\div}_{\g_2} \Y+{\d-1\over 2}\,\Y^2=0 \cr
       \end{cases}
       \Rightarrow
       {\div}_{\g_1}(\X+\Y)+{\d-1\over 2}\,(\X+\Y)^2=0\,.
                      \end{equation*}
From equation~\eqref{differenceoftwosingularoperators}, follows  a `cocycle condition': that
the sum of the left-hand side of first two equations is equal to the left-hand side of the third equation.
\end{remark}

\begin{remark}\label{interestingoperatoronfucntionsex}
Let $\rh$ be an arbitrary non-vanishing volume form.
Using the operator $\Delta([\g])=\Delta_\sing(\SS,\g)$ one can consider the second order operator $\Delta$ on functions,
\begin{equation*}
\Delta f=
 \rh^{-{1+\d\over 2}}\Delta\left([\g]\right)\left(\rh^{1-\d\over 2}f(x)\right)\,,
 \end{equation*}
depending on a volume element $\rh$.
A calculation gives an explicit formula:
                   \begin{equation*}
                     \Delta f={1\over 2}\Bigl(
                         S^{ab}\p_a\p_b+
                \p_bS^{ba}\p_a+
                (\d-1)\gamma_\rh^a\p_a+
                U_{\SS,\g}-U_{\SS,\g^\rh}
             \Bigr)f
              \end{equation*}
Here $\g_\rh\colon \,\gamma_a=-\p_a\log \rho$ is the flat connection  defined by the volume element $\rh=\rho(x)\D{x}$,
  $\gamma^a=S^{ab}\gamma_b$ and $U_{\SS,\g}$ is the ``pseudoscalar" part \eqref{pseudoscalarpart} of the operator
  \eqref{operatorsingular2}.  The difference  $U_{\SS,\g}-U_{\SS,\g_\rh}$
  is a density of weight $\d$ (see Corollary \ref{corollary2} and
  equation \eqref{differenceoftwosingularoperators}).
\end{remark}

We shall consider now examples of groupoids $C_\SS$ and the corresponding operators $\Delta_\sing(\SS,\g)$.

\subsection {The groupoid $C_\SS$ for a Riemannian manifold}

Let $M$ be a   manifold equipped with a Riemannian metric $G$. (As always we suppose that $M$ is an orientable compact manifold with a chosen oriented atlas). The Riemannian metric defines a principal symbol $\SS=G^{-1}$. In local coordinates $S^{ab}=g^{ab}$, where $G=g_{ab}dx^adx^b$.  It is a principal symbol of an operator of weight $\d=0$.

Let $\g$ be an arbitrary connection on volume forms. The differential operator $\Delta=\Delta_{\sing} (G^{-1},\g)$ of weight $\d=0$
maps half-densities to half-densities, $\Delta_{\sing} (G^{-1},\g)\colon \F_{1\over 2}\to \F_{1\over 2}$.  According to equations \eqref{operatorsingular2} and \eqref{pseudoscalarpart},
this operator equals to
\begin{equation*}
                   \Delta_{\sing} (G^{-1},\g)={1\over 2}\left(
                   g^{ab}\p_a\p_b+\p_bg^{ba}\p_a+{1\over 2}\,\p_a \gamma^a-{1\over 4}\,\gamma_a\gamma^a
                   \right).
\end{equation*}
We come  to the groupoid
\begin{equation*}
 C_\SS=\left\{{\g{\buildrel \X\over \longrightarrow}
 \,\g^\pr} \ \Bigl|\Bigr. \  {\div}_{\g} \X-{1\over 2}\,\X^2=0 \right\}\,
\end{equation*}
and the operator on half-densities depending on a class of connections
 \begin{equation*}
      \Delta\left([\g]\right)={1\over 2}\left(
                   g^{ab}\p_a\p_b+\p_bg^{ba}\p_a+{1\over 2}\,\p_a \gamma^a-{1\over 4}\,\gamma_a\gamma^a
                   \right).
\end{equation*}

On a Riemannian manifold one can consider the distinguished Levi-Civita connection. This connection
defines the connection   $\g^G$ on volume forms,
such that $\gamma^G_a=-\Gamma_{ab}^b=-\p_a\log \sqrt{\det g}$, where
$\Gamma^a_{bc}$ is the Christoffel symbol for the  Levi-Civita connection.
(We refer to  this connection on volume forms, also as to the Levi-Civita connection.)
Consider the orbit $[\g^G]$ in the groupoid
$C_G$ of the Levi-Civita connection $\g^G$.
This orbit defines a distinguished operator
on half-densities on a Riemannian manifold:
             \begin{equation*}\label{distinguishedoperator}
                        \Delta=\Delta_G\left([\g^G]\right)\,.
                              \end{equation*}

One can always choose special  local coordinates $(x^a)$ such that in these coordinates
$\det g=1$. In these local coordinates $\gamma^G_a=0$ and
the distinguished operator $\Delta$ on half-densities has the form
                        \begin{equation*}\label{distinguishedoperatorlocal}
                        \Delta={1\over 2}\bigl(
                   g^{ab}\p_a\p_b+\p_bg^{ba}\p_a\bigr)\,,
                              \end{equation*}
i.e., on a half-density $\ss=s(x)\D{x}^{1\over 2}$,
\begin{equation*}
    \Delta \ss= {1\over 2}\Bigl(
                   \p_b\left(g^{ba}\p_a s\right)\Bigr)\D{x}^{1\over 2}\,.
\end{equation*}
The differential equation
               $$
   {\div}_{\g} \X-{1\over 2}\,\X^2=0
             $$
defining the groupoid $C_G$ has the following form in these coordinates:
                   \begin{equation*}\label{interestingequation}
                   {\p X^a(x)\over \p x^a}-{1\over 2}\,X^a(x)X_a(x)=0\,.
                   \end{equation*}
All connections $\g$ of the form  $\gamma_a(x)=X_a(x)$ in these special coordinates, where $X_a(x)$ is a solution of this differential equation, belong to the orbit $[\g^G]$.

The operator $\Delta\left([\g^G]\right)$ belongs, in particular, to the canonical pencil associated with the Laplace--Beltrami operator
(see Example~\ref{beltrami-laplaceoperatorondensities2}).

 Let,  on the other hand,   $\g$ be an arbitrary connection and let $\rh$ be an arbitrary non-vanishing volume form on
a Riemannian manifold $M$. One can assign to a volume form $\rh$ the flat connection $\g^\rh\colon \gamma^\rh_a=-\p_a\log \rho$.
Consider the operator ${1\over \sqrt \rh }\Delta\left([\g]\right)\sqrt \rh$ on functions (see Remark \ref{interestingoperatoronfucntionsex}.)   We come to an scalar operator on functions
            \begin{equation*}
         \Delta f={1\over 2}\left(\p_a\left(g^{ab}\p_b f\right)-\gamma^{\rh a}\p_a+R\right)
                  \end{equation*}
   where the scalar function $R$ equals to
                     $$
              R=U_{G,\g}-U_{G,\g^\rh}={1\over 2}\,\div \X-{1\over 4}\X^2\,.
                     $$
   Here the vector field $\X$ is defined as the difference of the connections: $\X=\g-\g^\rh$, and $\div\X$ denotes the Riemannian divergence.

It will be interesting to compare the formulae of this subsection  with the constructions
in paper~\cite{BatBer2} for a case of Riemannian structure.

\section{Further examples}

Our next example is a groupoid of connections arising on an odd symplectic supermanifold.
Before discussing it, let us discuss briefly what happens  if we consider supermanifolds instead of manifolds.

\subsection {The supermanifold case}

 Let $M$ be an $n|m$-dimensional supermanifold. Denote local coordinates
   on $M$ by $z^A=(x^a,\theta^\a)$ ($a=1,\dots,n\,; \a=1,\dots,m$).
   Here $x^a$ are even coordinates and $\theta^\a$ are odd coordinates:
                       $z^Az^B=(-1)^{p(A)p(B)}z^Bz^A$,
   where $p(z^A)$, or shortly $p(A)$ is the parity of a coordinate $z^A$, so $p(x^a)=0$ and $p(\theta^\a)=1$.

   We would like to study second order linear differential operators
   $\Delta=S^{AB}\p_A\p_B+\dots$. The principal symbol of this operator is a
   supersymmetric contravariant tensor field   $\SS=(S^{AB})$. This field may be even or
   odd:
                      $$
        S^{AB}=(-1)^{p(A)p(B)}S^{BA},\,\, p(S^{AB})=p(\SS)+p(A)+p(B)\,.
                      $$
The analysis of second order operators can be performed in the supercase
in a way similar to the usual case. We  only have to be careful with signs. E.g., equation~\eqref{canonicalpencil} for a canonical pencil  has to be rewritten in the following way:
\begin{multline}\label{canonicalpencilinsuper}
\Delta_\l=\frac{t^\delta}{2}\Bigl(
    S^{AB}(x)\p_B\p_A+(-1)^{p(A)p(\SS+1)}\p_BS^{BA}\p_A \ +\\
                             +
    \left(2\l+\delta-1\right) \gamma^A \p_A+(-1)^{p(A)p(\SS+1)}\l\p_A\gamma^A +
         \l \left(\l+\delta-1\right)\theta
             \Bigr)\,.
\end{multline}
Here $\Delta$ is an even (odd) operator if the principal symbol $\SS$ is an even (odd) tensor density
(see~\cite{KhVor2} for details).

If $\SS$ is an even tensor field and it is non-degenerate, then it defines a Riemannian structure
on the (super)manifold $M$.   We come to the groupoid  $C_\SS$ in the  same way as in
the case of ordinary Riemannian manifolds considered in the previous section.
(We only have to take care of signs arising in calculations.)
In particular, for an even Riemannian supermanifold there
exists a distinguished   connection (the Levi-Civita connection), which canonically induces a unique connection on volume forms.
This  is a flat connection corresponding to the canonical volume form:
\begin{equation}\label{volumeformforsuperandlevicivitaconnection}
 \rh_g=\sqrt {\Ber(g_{AB})}\,|D(z)|\,,\quad
 \gamma_A=-\p_A\log \rho(z)=-(-1)^B\Gamma^B_{BA}\,.
\end{equation}
 Here $g_{AB}$ is the covariant tensor defining the Riemannian structure,
 ($S^{AB}=g^{AB}$) and $\Gamma^A_{BC}$ is the Christoffel symbol  of the Levi-Civita connection
for this Riemannian structure. $\Ber (g_{AB})$ is  the Berezinian of the matrix $g_{AB}$. It is the super analog of determinant. The matrix $g_{AB}$
 is an $n|m\times n|m$ even matrix and its Berezinian is given by the formula
          \begin{equation}\label{Berezinian}
           \Ber(g_{AB})=
              \Ber
              \begin{pmatrix}
              g_{ab} & g_{a\beta}\cr
              g_{\a b} &g_{\a\beta}
              \end{pmatrix}=
              \frac{\det \left(g_{ab}-g_{a\gamma}g^{\gamma\delta}g_{\delta b}\right)}{\det (g_{\a\beta})}\,.
              \end{equation}
(Here, as usual, $g^{\gamma\delta}$ stands  for the matrix inverse to the matrix $g_{\gamma\delta}$.)

The situation is essentially different   if $\SS=(S^{AB})$ is an odd
tensor field and, respectively,
$\Delta=S^{AB}\p_A\p_B+\dots$ is an odd operator.
In this case one comes naturally to an odd Poisson structure on the supermanifold $M$, provided
the tensor $\SS$ obeys additional constraints.

Namely, consider the cotangent bundle $T^*M$ to the supermanifold $M$ with local coordinates
$(z^A,p_A)$, where $p_A$ are the fiber coordinates  dual to coordinates $z^A$ (the variables $p_A$ transform as the partial derivatives $\lder{}{z^A}$).
A supersymmetric contravariant tensor field $\SS=S^{AB}$ defines a quadratic master-Hamiltonian, an odd function
  $H_\SS={1\over 2}S^{AB}p_Ap_B$ on the cotangent bundle $T^*M$.
This quadratic master-Hamiltonian defines an odd bracket on functions on $M$
as a derived bracket:
\begin{equation}\label{derivedbracket1}
                     \{f,g\}=\left(\left(f,H_\SS\right),g\right)),\qquad
                     p\left(\{f,g\}\right)=p(f)+p(g)+1\,.
\end{equation}
Here $(\,,\,)$ is the canonical Poisson bracket on the cotangent bundle $T^*M$.
The odd derived bracket is antisymmetric with respect to  shifted parity  and it obeys the Leibniz rule:
              $$
      \{f,g\}=-(-1)^{(p(f)+1)(p(g)+1)}\{g,f\},\quad \{f,gh\}=\{f,g\}h+(-1)^{(p(f) +1)p(g)}g\{f,h\}\,.
              $$
This odd derived bracket is  an odd Poisson bracket if it obeys the Jacobi identity
                     \begin{equation}\label{jacobiidentity}
(-1)^{(p(f) +1)(p(h) +1)}\{\{f ,g \},h \}+(-1)^{(p(g) +1)(p(f) +1)}\{\{g ,h \},f \}
+(-1)^{(p(h) +1)(p(g) +1)}\{\{h , f\},g \}=0\,.
                      \end{equation}
It is a beautiful fact that the condition that the derived bracket \eqref{derivedbracket1} obeys the Jacobi identity
can be formulated as a quadratic condition $\left(H,H\right)=0$ for the master-Hamiltonian:
                \begin{equation}\label{derivedbracket2}
                     \left(H_\SS,H_\SS\right)=0\,
                      \Leftrightarrow
                     \text{Jacobi identity
                     for the derived bracket $\{\,,\,\}$.}
                      \end{equation}
(See~\cite{KhVor1} for details.\footnote{If $\SS$ is an even tensor field--- the case of Riemannian geometry---the master-Hamiltonian
 $H$ is an even function and the Jacobi identity is trivial, see \cite{KhVor1} and
\cite{KhVor2} for details.})

From now on  suppose that the tensor field $\SS$ is odd and it defines an odd Poisson bracket
  on the supermanifold $M$,
  i.e., that the relation~\eqref{derivedbracket2} holds.
  This odd Poisson bracket corresponds to an odd symplectic structure
  in the case if the bracket is non-degenerate, i.e., the odd tensor   $\SS$ is   non-degenerate.
  The condition of non-degeneracy  means that there exists the inverse covariant tensor $S_{BC}$:
  $S^{AB}S_{BC}=\delta^B_C$.  Since the matrix $S^{AB}$ is an odd matrix
  ($p(\SS^{AB})=p(A)+p(B)+1$),
  this implies that the matrix $S^{AB}$ has equal  even and odd dimensions.
  We come to the conclusion that for an odd symplectic supermanifold the even and odd dimensions have to coincide.
  It is necessarily $n|n$-dimensional.

The basic example of an odd symplectic supermanifold is the following.
     for an arbitrary ordinary manifold $M$ consider its cotangent bundle $T^*M$ and change parity of the fibres
     in this bundle.  We arrive at the odd symplectic supermanifold $\Pi T^*M$.
With arbitrary local coordinates  $x^a$ on
     $M$, one can associate local coordinates $(x^a,\theta_a)$ on $\Pi T^*M$,
     where the odd coordinates $\theta_a$ transform in the same way as $\p_a$:
\begin{equation}\label{trasnformationoflocalcoordinatesinsuperpsace}
                       x^{a^\pr}=x^{a^\pr}(x)\,,\qquad
               \theta_{a^\pr}=\der{x^a}{x^{a^\pr}}\,\theta_a\,.
\end{equation}
  In these local coordinates the canonical non-degenerate odd Poisson bracket is  defined by the relations
                 \begin{equation}\label{Darbouxcoordinates}
                 \{x^a,x^b\}=0\,, \quad \{x^a,\theta_b\}=\delta^a_b\,, \quad \{\theta_a,\theta_b\}=0\,.
                  \end{equation}
(These relations are invariant with respect to coordinate
transformations \eqref{trasnformationoflocalcoordinatesinsuperpsace}.)

\begin{remark}\label{symplecctomorphic}
Every odd symplectic supermanifold $E$ is symplectomorphic to the anticotangent bundle $\Pi T^*M$ of an ordinary manifold $M$ (i.e., the cotangent bundle with reversed parity in the fibers).
One may take as $M$  a purely even Lagrangian surface in $E$. (See~\cite{Kh3} for details.)
Therefore,  one  can always consider  $\Pi T^*M$ instead of an odd symplectic supermanifold $E$. Note that $\Pi T^*M$ possesses an extra structure, viz., that of a vector bundle over $M$. A practical
difference between  $E$ and $\Pi T^*M$ is that for an `abstract'  odd symplectic supermanifold $E$ one may consider arbitrary parity preserving
coordinate transformations of local coordinates $x$ and $\theta$, while for $\Pi T^*M$ there is a privileged class of
transformations of the form~\eqref{trasnformationoflocalcoordinatesinsuperpsace}, typical for a vector bundle.
\end{remark}

\subsection {The groupoid $C_\SS$ for an odd symplectic supermanifold}

Let $E$ be an arbitrary $(n|n)$-dimensional odd symplectic supermanifold, where the
odd symplectic structure and respectively the odd non-degenerate Poisson structure are defined by
a non-degenerate contravariant supersymmetric  odd tensor field $\SS=(S^{AB})$.\footnote{The following construction of groupoid
   is obviously valid in the general Poisson case, but in this subsection we are mainly interested
   in the odd symplectic case}.   We shall study odd second order   operators
   $\Delta={1\over 2}S^{AB}+\dots$ of weight $\d=0$.

Let $\g$ be an arbitrary connection on volume forms.
   The differential operator $\Delta=\Delta_{\sing} (\SS,\g)$ of weight $\d=0$
   with the principal symbol $\SS$ is defined by equation \eqref{operatorsingular2}.
It transforms half-densities to half-densities, $\Delta_{\sing} (\SS,\g)\colon \F_{1\over 2}\to \F_{1\over 2}$.
From equations~\eqref{operatorsingular2}, \eqref{pseudoscalarpart}
   and \eqref{canonicalpencilinsuper}, we obtain
\begin{equation}\label{bvsingularoperator}
                   \Delta_{\sing} (\SS,\g)=\frac{1}{2}\left(S^{AB}\p_B\p_A+\p_Bg^{BA}\p_a+\frac{1}{2}\,\p_A \gamma^A-\frac{1}{4}\,\gamma_A\gamma^A\right).
\end{equation}
We come  to the groupoid
                         \begin{equation*}
                         C_\SS=\left\{{\g\,{\buildrel \X\over \longrightarrow}\, \g^\pr}
 \ \Bigl|\Bigr. \ {\div}_{\g} \X-{1\over 2}\,\X^2=0
                            \right\}
                          \end{equation*}
and to the operator on half-densities depending on a  class of connections (a groupoid orbit)
                 \begin{equation}\label{bvsingonclasses}
      \Delta\left([\g]\right)={1\over 2}\left(
                   S^{AB}\p_B\p_A+\p_Bg^{BA}\p_a+U_\SS\left([\g]\right)\right)
                   \,,\,\,\text{where }\,  U_\SS\left([\g]\right)=
                   \frac{1}{2}\,\p_A \gamma^A-\frac{1}{4}\,\gamma_A\gamma^A\,.
                 \end{equation}
It is here where the similarity with the Riemannian case finishes.
On a Riemannian manifold, one can consider the canonical volume element and the
distinguished Levi-Civita connection, which give
the canonical flat connection $\g$ on volume forms (see equation \eqref{volumeformforsuperandlevicivitaconnection}).
On an odd symplectic supermanifold, there is no
canonical volume element\footnote {A naive generalization of formulae
\eqref{volumeformforsuperandlevicivitaconnection} and \eqref{Berezinian}
is not possible because, in particular, $S^{AB}$ is not an even matrix.} and there is no
distinguished affine connection.

On the other hand, it turns out that for an odd symplectic supermanifold one can construct
a distinguished \emph{equivalence class of connections} on volume forms (a distinguished orbit of the groupoid $C_\SS$).
Namely, consider the equation
              \begin{equation}\label{bvequation2}
   {\div}_{\g} \X-\frac{1}{2}\,\X^2=0\,,
               \end{equation}
which defines the groupoid $C_\SS$.
According to equations \eqref{differenceoftwosingularoperators}, \eqref{bvsingularoperator}
and \eqref{bvsingonclasses}, for the operators   $\Delta\left([\g]\right)$
acting on half-densities, we have
                   \begin{equation}\label{bv2}
                   \Delta\left([\g']\right)-\Delta\left([\g]\right)=
                   \Delta_{\sing} (\SS,\g')-\Delta_{\sing} (\SS,\g)=
                   {1\over 4}\left({\div}_{\g} \X-\frac{1}{2}\,\X^2\right)\,.
                   \end{equation}
We call  equation \eqref{bvequation2},  the \emph{Batalin-Vilkovisky equation}. Let us study it.

It is convenient to work in  Darboux coordinates.  Recall that local coordinates
  $z^A=(x^a,\theta_a)$ on an odd symplectic supermanifold $E$ are called \emph{Darboux coordinates}
  if  the pairwise odd Poisson brackets have the form~\eqref{Darbouxcoordinates}.
We say that a connection $\g$ is \emph{Darboux flat} if it vanishes in some Darboux coordinates.

\begin{lemma}\label{lemma}
Let $\g,\g'$ be two connections   that are both  Darboux flat. Then the arrow
$\grarrow$  belongs to the groupoid $C_\SS$, i.e., the Batalin-Vilkovisky equation
${\div}_{\g} \X-\frac{1}{2}\,\X^2=0$ holds for the covector field $\X=\g'-\g$.
\end{lemma}
(We shall prove this lemma later.)

{\small
\begin{remark}\label{localconsiderations}
In fact, lemma implies that the class of local   Darboux flat connections
defines a global  pseudoscalar function  $U_\SS$ in \eqref{bvsingularoperator}.
 Let $\{z^A_{(\a)}\}$ be an arbitrary
atlas of Darboux coordinates on $E$.
We say  that a collection
of local connections $\{\g_{(\a)}\}$ is adapted to the Darboux atlas
$\{z^A_{(\a)}\}$
if every local connection $\g_{(a)}$ (defined  in the chart $z^A_{(\a)}$) vanishes in these
local Darboux coordinates
$z^A_{(\a)}$.
  Let  $\{\g_{(\a)}\}$ and $\{\g'_{(\a')}\}$
be two families of  local connections adapted to Darboux atlases  $\{z^A_{(\a)}\}$
and $\{z^{A^\pr}_{(\a^\pr)}\}$ respectively. Then, due to Lemma, all the arrows
 ${\g_{(\a)}\buildrel \X\over \longrightarrow}\g_{(\a')}$
 $ {\g^\pr_{(\a)}\buildrel \X\over \longrightarrow}\, \g^\pr_{(\a')}$
 and
 ${\g_{(\a)}\buildrel \X\over \longrightarrow}\, \g^\pr_{(\a')}$
  belong to a local groupoid $C_\SS$ (if the charts
  $(z^A_{(\a)})$, $(z^A_{(\a')})$,
  $(z^{A'}_{(\a)})$ and $(z^{A'}_{(\a')})$ intersect).
This means that in spite of the fact that the family $\{\gamma_\a\}$
does not define a global connection, still
equations \eqref{bvequation2} hold locally
and the operator $\Delta=\Delta(\SS,\g_\a)$ globally exists.
(These considerations for a locally defined groupoid
can be performed for an arbitrary case.  One can consider the family of locally defined
connections $\{\g_a\}$ such that they define a global operator
\eqref{operatorsingular2}.)
On the other hand in the case of an odd symplectic supermanifold
there exists a global Darboux flat connection,
i.e., a connection
  $\g$ such in a vicinity of an arbitrary point this connection vanishes in some Darboux coordinates.
Let us show that.

Without loss of generality suppose that $E=\Pi T^*M$ (see Remark~\eqref{symplecctomorphic}.)
Let $\s$ be an arbitrary non-vanishing volume form on $M$  (we suppose that $M$ is orientable).
Choose an atlas $\{x^a_{(\a)}\}$  of local coordinates on $M$
such that $\s$ is the coordinate volume form in each chart,
i.e., $\sigma=dx^1_{(\a)}\wedge \dots dx^n_{(\a)}$.   Then consider the associated atlas
$\{x^a_{(\a)},\theta_a{(\a)}\}$ on the supermanifold $\Pi T^*M$, which is a Darboux atlas.
For this atlas, as well  as for the atlas
$\{x^a_{(\a)}\}$, the Jacobians of the coordinate transformations are equal to $1$.
Thus we have constructed a `special' Darboux atlas for $\Pi T^*M$  for which all the Jacobians
of the coordinate transformations are equal to $1$.  Therefore the coordinate volume element
$\rho=D(x,\theta)$ is globally defined. The components of the flat connection corresponding to this volume element
vanish in each coordinate chart. Hence we have defined a global   Darboux flat connection.
\end{remark}

}

We have arrived at the following proposition.

\begin{proposition}\label{canonicalclassofconnections}
For an odd symplectic supermanifold there exists
a distinguished  orbit of the groupoid of connections $C_\SS$. It is the class $[\g]$
of an   arbitrary Darboux flat connection  $\g$.
\end{proposition}

We   call this canonically defined class of connections ``the  class of Darboux flat connections''.

For any connection belonging to this canonical class of connections, the pseudoscalar function  $U_{\SS}$ in equation~\eqref{bvsingonclasses}
 vanishes in arbitrary Darboux coordinates.\footnote{This function vanishes not only for globally defined Darboux flat connection
 but for a  family of connections adapted to an arbitrary Darboux atlas (see Remark \ref{localconsiderations})}
 The operator  $\Delta=\Delta[\g]$ on half-densities corresponding to this class of connections
 has the following appearance in arbitrary Darboux
 coordinates $z^A=(x^a,\theta_b)$:
\begin{equation}\label{khudian}
                  \Delta={\p^2\over \p x^a\p \theta_a}\,.
\end{equation}
(This canonical operator on half-densities was introduced in \cite{Kh3}.)

Let us now prove Lemma~\ref{lemma}.

For an arbitrary non-vanishing volume form $\rh$ consider the operator
\begin{equation}\label{bvkhudian}
\Delta_\rh f=\frac{1}{2}\,{\div}_\rh \grad f\,.
\end{equation}
Here $\grad f$ is the Hamiltonian vector field $\{f,z^A\}{\p \over \p z^A}$ corresponding to a function $f$.
(Compare with \eqref{firstequation}.)
This is the famous Batalin--Vilkovisky odd Laplacian on functions.
If $z^A=(x^a,\theta_a)$ are Darboux coordinates and the volume form $\rh$
is the coordinate volume form, i.e.,  $\rh=D(x,\theta)$,
then the odd Laplacian  in these Darboux coordinates
has the appearance
           \begin{equation}\label{bvinitial}
           \Delta=\dder{}{x^a}{\theta_a}\,.
            \end{equation}
(This is the original form of the Batalin--Vilkovisky operator in~\cite{BatVyl1}.
For the geometric meaning of the BV operator  and for the explanation how  formulae~\eqref{bvkhudian} and \eqref{bvinitial} are related with the canonical operator~\eqref{khudian} on semidensities see~\cite{Kh1, Sch, Kh3}.)

Equation~\eqref{bvequation2} (the Batalin--Vilkovisky equation) characterizing the groupoid $C_{\SS}$ is related with the Batalin--Vilkovisky operator by the following identity:
\begin{equation}\label{identitysimple}
            -e^{\frac{F}{2}}  \Delta_\rh e^{-\frac{F}{2}}=
           \frac{1}{4}\,\left({\div}_\g \X-\frac{1}{2}\,\X^2\right)\,,
\end{equation}
where the  connection $\g$ is the flat connection induced by a volume form
($\gamma_a=-\p_a\log \rho$ ) and the vector field $\X$
is the Hamiltonian vector field of a function $F$.

We use this identity to prove the Lemma. Let the connection $\g$ vanish  in Darboux coordinates $z^A=(x^a,\theta_a)$ and the connection $\g^\pr$ vanishes in Darboux coordinates $z^{A'}=(x^{a'},\theta_{a'})$. Then  (compare with equation \eqref{transformationofgenuineconnection})
the connection $\g'$ has in Darboux coordinates $z^A=(x^a,\theta_a)$ the following form:
                    $$
\gamma^\pr_A=\der{z^{A^\pr}}{z^A}\Bigl(\gamma_{A^\pr}+\p_{A'}\log J\Bigr)\,,
                    $$
where $J$ is the Jacobian of transformation of the Darboux coordinates: $J=\Ber\der{(x,\theta)}{(x^\pr,\theta')}$  (see also formula \eqref{Berezinian}).
Hence for the arrow $\grarrow$ the components of the covector field $\X$ are
 $X_A=-\p_A\log \Ber\der{(x^\pr,\theta^\pr)}{(x,\theta)}$. Now we apply the identity~\eqref{identitysimple}, where $\gamma_a=0$, $\rh=D(x,\theta)$ is the  coordinate volume form, and the function $F$ is $-\log \Ber\der{(x^\pr,\theta^\pr)}{(x,\theta)}$. Using~\eqref{bvkhudian} and \eqref{bvinitial}, we arrive at
 \begin{equation*}
     \frac{1}{4}\,\left({\div}_\g \X-\frac{1}{2}\,\X^2\right)= -e^{\frac{F}{2}}\Delta_\rh e^{-\frac{F}{2}}=
                       -           \left(\sqrt {\Ber\der{(x,\theta)}{(x^\pr,\theta^\pr)}}\,\right)
             \dder{}{x^a}{\theta_a}
             \left(\sqrt{\Ber\der{(x^\pr,\theta^\pr)}{(x,\theta)}}\,\right)=0\,.
 \end{equation*}
The last equality  is the fundamental Batalin--Vilkovisky identity \cite{BatVyl2} (see also~\cite{KhVor1}), which is in the
core of the geometry of Batalin--Vilkovisky operator.

Lemma~\ref{lemma} is now proved.

\begin{remark}   The canonical operator~\eqref{khudian} allows to assign to every invertible half-density $\ss$
and  every invertible volume form $\rh$ two functions $\sigma_\ss$ and $\s_\rh$\,:
                      $$
                \sigma(\ss)={\Delta \ss\over \ss},\qquad \sigma (\rh)={\Delta\sqrt \rh\over \sqrt{\rh}}
                      $$
(see \cite{Kh3}).
In articles \cite{BatBer1, BatBer2},  Batalin and Bering were studying geometric
properties of the canonical operator
\eqref{khudian}.  In these considerations they were using the formula for
this operator in arbitrary coordinates suggested by Bering in~\cite{Ber}.
Trying to clarify the geometric meaning of this formula
and whilst analyzing the  meaning of the scalar function $\s(\rh)$,  they came to the following beautiful result:
if $\nabla$ is an arbitrary torsion-free affine connection on an odd symplectic supermanifold compatible with a volume element $\rh$, then
the scalar curvature of this connection is equal, up to a factor,  to the function $\sigma(\rh)$.
\end{remark}

\subsection{The groupoid $C_\SS$ for the line}
We return here to the simplest possible manifold---the real line.
A  tensor density $\SS$ of rank $2$ and  weight $\delta$ on the real line $\RR$ can be identified with a density of   weight $\delta-2$:
\begin{equation*}
    \SS=S(x)\D{x}^\d\,\p_x\otimes \p_x\sim S(x)\D{x}^{\d-2}\,.
\end{equation*}
Consider on $\RR$ the  vector density $\D{x}\p_x$ ,which is invariant with respect to change of coordinates.
Its square defines the invariant tensor density  $\SS_\RR=\D{x}^2\p_x\otimes \p_x$ of weight $\d=2$. (In the above identification, both correspond to the function $1$.)

We see that on the line there are canonical pencils of second order operators
of  weight $\d=2$\,:
\begin{equation*}
    \Delta_\l=\D{x}^2\Bigl(\p_x^2+\dots\Bigr)
\end{equation*}
with the distinguished principal symbol $\SS_\RR$.
The operator $\Delta(\g):=\Delta_\sing(\SS_\RR,\g)$ belonging to this pencil acts on  densities of weight ${1-\d\over 2}=-{1\over 2}$ and transforms them into densities of weight ${1+\d\over 2}={3\over 2}$. According to~\eqref{operatorsingular2}, it has the following form:
                      $$
\Delta(\g)\colon\, \Psi(x)|D(x)|^{-{1\over 2}}\mapsto
              \Phi(x)\D{x}^{3\over 2}=
              {1\over 2}\left({\p^2\Psi(x)\over \p x^2}+U(x)\Psi(x)\right)|D(x)|^{{3\over 2}}\,,
                   $$
where,  according to equation \eqref{pseudoscalarpart},
                \begin{equation}\label{sturmoperator}
             U_\g(x)=-{1\over 4}\left(\gamma_x+{1\over 2}\,\gamma^2\right)\D{x}^2.
                 \end{equation}
The operator $\Delta(\g)$ is the Sturm--Liouville operator, well known to experts in projective geometry
and integrable systems (see, e.g., \cite{HitchSeg} or \cite{OvsTab}).\footnote
{An operator $\Delta$ corresponds to the curve $t\mapsto [u_1(t):u_2(t)]$,
$\RR\ \to \RR P^1$ in projective line defined by two solutions of equation $\Delta u=0$.}

We see that in this case the difference of the operators is
\begin{equation*}
    \Delta(\g')-\Delta(\g)=
            -{1\over 4}\left(\gamma'_x+{1\over 2}\, {\gamma^{\pr}}^2\right)\D{x}^2+
            {1\over 4}\left(\gamma_x+{1\over 2} \,\gamma^2\right)\D{x}^2=\\
            -{1\over 4}\left({\div}_\g\X+\frac{1}{2}\,\X^2\right).
\end{equation*}
Here $\X=(\gamma^\pr-\gamma)\D{x}^2\p_x$ is a vector density of   weight $\d=2$, which can be identified with a density of weight $1$ on the line.
(compare with formulae\eqref{differenceoftwosingularoperators} and \eqref{bv2}).

By applying the general formulae \eqref{definitionofgroupoid} to the case of the distinguished principal symbol $\SS_\RR$ on the line, we arrive at the following canonical groupoid $C_\RR$ for the line:
\begin{equation*}\label{groupidonline}
   C_\RR= \left\{{\g{\buildrel \X\over \longrightarrow}
 \,\g^\pr} \ \Bigl|\Bigr. \  \Delta(\g')=\Delta(\g), \ \text{i.e.,} \   U_{\g^\pr}=U_\g \right\}=
 \left\{{\g{\buildrel \X\over \longrightarrow}
 \,\g^\pr} \ \Bigl|\Bigr. \  {\div}_\g\X+\frac{1}{2}\,\X^2=0\right\}\,,
\end{equation*}
where $\Delta(\g)$ is the Sturm-Liouville operator~\eqref{sturmoperator}.
The Sturm-Liouville operator depends on the orbit $[\g]$ of a connection $\g$ in the groupoid $C_\RR$.


Let us analyze the equation ${\div}_\g\X+\frac{1}{2}$ defining the canonical groupoid $C_\RR$
and compare it with a cocycle related with the operator.
For the covector field  $\g^\pr-\g=a(x)dx$, the corresponding vector density equals to
$\SS_\RR\bigl(a(x) dx\bigr)=a(x)\D{x}^2\p_x$. Hence $\X^2=a^2(x)\D{x}^2$ and
 ${\div}_\g\X=(a_x+\gamma a)\D{x}^2$. We arrive at the equation
\begin{equation*}
           {\div}_\g\X+\frac{1}{2}\,\X^2=
            \left(a_x+\gamma a+\frac{1}{2}\,a^2\right)\D{x}^2=0\,.
\end{equation*}
To solve this differential equation, choose a coordinate on $\RR$ such that the coefficient $\gamma$ vanishes in this coordinate.
Then
               \begin{equation}\label{solutionofdifferentialequation}
  \X={2dx\over C+x}, \qquad \hbox {where $C$ is a constant}\,.
               \end{equation}
On the other hand,  analyze the action of diffeomorphisms on the connection $\g$ and on the
Sturm-Liouville operator~\eqref{sturmoperator}.  Let $f=f(x)$ be a diffeomorphism of $\RR$.
(We consider the compactified $\RR$, i.e., $S^1=\RR P^1$ and the diffeomorphisms preserving orientation.)
The new connection $\g^{(f)}$ is equal  to
$y_x\left(\g\big\vert_{y(x)}+\left(\log x_y\right)_x\right)dx$
and the covector field $\g^{(f)}-\g$ equals
                     $$
                     \X^{(f)}=\g^{(f)}-\g=\gamma(y(x))dy+\left(\log x_y\right)_y dy-\gamma(x)dx\,.
                      $$
We come to a cocycle  on the group of diffeomorphisms:
\begin{equation}\label{cocycle}
    c_\g(f)=\Delta(\g^f)-\Delta(\g)=
                {1\over 4}\left(U_{\g^f}-U_{\g}\right)=\\
                -{1\over 4}\left(\div\X^{(f)}+{1\over 2}\left(\X^{(f)}\right)^2\right)\,.
\end{equation}
In a coordinate such that $\gamma=0$, we have $\X^{(f)}=\left(\log x_y\right)_y dy$. Combining this
with the solution \eqref{solutionofdifferentialequation}, we arrive at the
equation $\left(\log x_y\right)_x={2\over C+x}$. Solving  this equation gives
         $$
   \div\X^{(f)}+{1\over 2}\left(\X^{(f)}\right)^2=0\Leftrightarrow  y={ax+b\over cx+d}
    \text{(a projective transformation)}\,.
         $$
The cocycle~\eqref{cocycle} is the coboundary on the space of second order operators and it is a non-trivial cocycle
on the space of densities of   weight $2$. This cocycle vanishes on projective transformations.
This is a well-known cocycle related with the Schwarzian derivative
(see  book \cite{OvsTab} and citations therein):
\begin{multline*}
     c_\g(f)=\Delta(\g^f)-\Delta(\g)=\frac{1}{2}\left(U_{\g^f}-U_{\g}\right)=
      -\frac{1}{4}\left(\div \X^{(f)}+\frac{1}{2}\,\left(\X^{(f)}\right)^2\right)=\\
             -\frac{1}{4}\Bigl(U_\g(y)|D(y)|^2+{\cal S}\left[x(y)\right]|D(y)|^2-U_\g(x)\D{x}^2\Bigr)\,.
\end{multline*}
Here
 \begin{equation*}
   {\cal S}[x(y)]={x_{yyy}\over x_y}-{3\over 2}\left({x_{yy}\over x_y}\right)^2
\end{equation*}
is the \emph{Schwarzian} (or \emph{Schwarzian derivative}) of the transformation $x=x(y)$.
If $\gamma=0$ in the coordinate $x$, then $c(f)=-\frac{1}{4}\,{\cal S}[x(y)]|D(y)|^2$.

\subsection{The Lie algebroids of the Batalin-Vilkovisky groupoid and the groupoid  $C_\SS$}
  In work \cite{KhVor1}, we considered in particular the following ``Batalin-Vilkovisky
   groupoid'' on an odd Poisson manifold: the objects
are (non-vanishing) volume forms and
the arrows $\rh{\buildrel J\over \longrightarrow }\rh^\pr$, where $J={\rh^\pr\over \rh}$,
are specified by
the Batalin-Vilkovisky equation $\Delta_\rh \sqrt J=0$.
The operator $\Delta_\rh$ is defined by formula~\eqref{bvkhudian}.
Assign to each arrow  $\rh{\buildrel J\over \longrightarrow }\rh^\pr$
the arrow  $\grarrow$  of the groupoid $C_\SS$ such that the connections $\g,\g^\pr$
are defined by the volume forms
$\rh,\rh^\pr$ respectively ($\gamma_a=-\p_a\log \rho$ and $\gamma'_a=-\p_a\log \rho'$).
Then it follows from equation \eqref{identitysimple}
that the Batalin--Vilkovisky groupoid is a subgroupoid of the groupoid $C_\SS$.

Both the Batalin--Vilkovisky groupoid
and the groupoid of connections $C_\SS$ considered here can be regarded as
Lie groupoids over infinite-dimensional manifolds, which are the
space $\Vol^{\times}(M)$ of the non-degenerate volume forms
and the space $\Conx(M)$ of the connections on volume forms on a
manifold $M$ respectively. (The space $\Conx(M)$ was denoted $\A$ above.)  The corresponding Lie algebroids can be
described as follows.

For the Batalin--Vilkovisky groupoid, the Lie algebroid is the vector
bundle over the (infinite-dimensional) manifold
$\Vol^{\times}(M)$ whose  the fiber over the point
$\boldsymbol{\rho}$ is the vector space of all solutions of the
equation
$\Delta_{\boldsymbol{\rho}} F=0$,
where $F\in C^{\infty}(M)$. The anchor is tautological: it sends a
function $F$ to the infinitesimal shift $\boldsymbol{\rho} \mapsto
\boldsymbol{\rho} +\vare \boldsymbol{\rho} F$. A section of this
bundle is a functional $F[\boldsymbol{\rho}]$ of a non-vanishing volume form with
values in functions on $M$, such that for each $\boldsymbol{\rho}=\rho(x)D(x)$,
the above equation is satisfied. The Lie bracket is the restriction
of the canonical commutator of vector fields on
$\Vol^{\times}(M)$ and can be expressed by the explicit
formula
$$
            [F,G] [\boldsymbol{\rho};x] = \int_M\!\! D(y)\, \rho(y) \left(
F[\boldsymbol{\rho};y]\frac{\delta G[\boldsymbol{\rho};x]}{\delta
\rho(y)}  -
G[\boldsymbol{\rho};y]\frac{\delta F[\boldsymbol{\rho};x]}{\delta
\rho(y)}\right)\,.
$$
Here we write $F[\boldsymbol{\rho};x]$ for the value of
$F[{\rh}] \in C^{\infty}(M)$ at $x\in M$.

For the groupoid of connections $C_{\SS}$ (for a fixed tensor density $S^{ab}$
of weight $\delta$), the Lie algebroid is the vector bundle over
$\Conx(M)$ whose fiber over $\boldsymbol{\gamma}\in
\Conx(M)$ is the vector space of all solutions of the
equation
${\rm div}_\g \X=0$. A section is a functional of
connections taking values in these vector spaces. The anchor is
tautological: it sends a covector field $X_a$ to the infinitesimal
shift of the connection: $\gamma_a\mapsto \gamma_a+\vare X_a$. The
Lie bracket can be expressed by the formula
$$
        [X,Y]_a[\boldsymbol{\gamma};x] = \int_M \!\! D(y) \left(
X_b[\boldsymbol{\gamma};y] \,\frac{\delta
Y_a[\boldsymbol{\gamma};x]}{\delta \gamma_b(y)}-
    Y_b[\boldsymbol{\gamma};y]\, \frac{\delta
X_a[\boldsymbol{\gamma};x]}{\delta \gamma_b(y)} \right)\,.
$$
(In the case of supermanifolds, the formulae for the brackets will contain extra
signs.) Note that, since the groupoids in question are subgroupoids of
trivial (or pair) groupoids,  these Lie algebroids are subalgebroids
of the respective tangent bundles.

\subsection {Invariant densities on  submanifolds of codimension $1|1$
 in an odd symplectic supermanifold and mean curvature}

In the previous examples we considered  second order operators depending on   equivalence classes of connections
on volume forms.
In particular, we considered the   class of the Darboux flat connections on an odd symplectic supermanifold
(see   Proposition \ref{canonicalclassofconnections}) and with the help  of this class redefined the canonical
operator \eqref{khudian}.

Here we shall consider another example of a geometric construction  that depends on second order derivatives and a class of connection, but now it will be an equivalence class of Darboux flat  affine connections rather than of connections on volume forms.

Let $E$ be an odd symplectic supermanifold equipped with a volume element $\rh$.
Let $C$ be a non-degenerate submanifold  of codimension $1|1$ in $E$
(meaning that induced pre-symplectic structure on $C$ is non-degenerate).
We call such a submanifold, a ``hypersurface".

For an arbitrary affine connection $\nabla$ on $E$ and an arbitrary vector field $\Psi$  consider the following object:
\begin{equation}\label{constr1}
     A(\nabla,\Psi)=\Tr\left(\Pi\left(\nabla \Psi\right)\right)-{\div}_\rho\Psi\,,
\end{equation}
where $\Pi$ is the projector on $1|1$-dimensional planes symplecto-orthogonal to the hypersurface $C$ at the points of this hypersurface.
We define these objects in a vicinity of $C$.  (Formula~\eqref{constr1} appeared in discussions with our student O.~Little.)
Let the vector field $\Psi$ be symplecto-orthogonal to the hypersurface $C$.  Then one can see that, at the points of $C$,
\begin{equation}\label{lastrelation}
      A(\nabla, f\Psi)=fA(\nabla, \Psi)
\end{equation}
for an arbitrary function $f$. Therefore  $A(\nabla,\Psi)$ is well-defined on $C$ even if the vector field $\Psi$  is defined only at $C$ provided $\Psi$ is symplecto-orthogonal to $C$. This object  is interesting because it is related with the canonical vector valued half-density and the canonical scalar half-density on the manifold $C$ (see~\cite{Kh2} for details.)

Namely, let $\Psi=\Psi^A\p_A$ be a vector field on the hypersurface $C$ symplecto-orthogonal to $C$. From now on we suppose that  it also obeys the following additional conditions
\begin{itemize}
 \item  $\Psi$ is  an odd vector field:  $p(\Psi)=p(\Psi^A)+p(A)=1$,
 \item  $\Psi$ is non-degenerate, i.e., at least one of its components is non-nilpotent,
  \item $\o(\Psi,\Psi)=0$, where $\o$ is the symplectic form in $E$.
\end{itemize}
One can see that these conditions uniquely define the vector field $\Psi$ at  every point of $C$ up to a factor (an invertible function).

Consider now the following volume form $\r_\Psi$ on $C$.
Let ${\bf H}$ be an even vector field on the hypersurface $C$ such that it is symplecto-orthogonal to $C$ and $\o({\bf H},\Psi)=1$.
Define a  half-density $\rh_\Psi$ on the  hypersurface $C$ by the condition that for an arbitrary
tangent frame  $\{\e_1,\dots \e_{n-1}; \f_1,\dots,\f_{n-1}\}$,
             $$
        \rh_\Psi\left(\e_1,\dots \e_{n-1}; \f_1,\dots,\f_{n-1}\right)=
        \rh \left(\e_1,\dots \e_{n-1}, {\bf H}; \f_1,\dots,\f_{n-1},\Psi\right)\,.
            $$
(Here $\e_1,\dots \e_{n-1}$ are even basis vectors and $\f_1,\dots,\f_{n-1}$ are odd basis vectors in a tangent space to $C$.)
Using formula \eqref{Berezinian} for Berezinian and relation
\eqref{lastrelation} one can see that for an arbitrary
function $f$,
$$
   \rh_{ f\psi}=\frac{1}{f^2}\, \rh_{\Psi}\,.
$$
We come to the conclusion that the odd vector-valued half-density $\Psi\sqrt{\rh_\Psi}$ is well-defined  on the hypersurface $C$. By applying equation \eqref{constr1}, we obtain a well-defined half-density on the hypersurface $C$:
                    $$
        \ss_{_C}( \nabla) =A(\nabla, \Psi)\sqrt {\rh_\Psi}\,.
                    $$
This half-density depends only on  an affine connection $\nabla$.

We say that an affine  symmetric connection $\nabla$ on $E$
with the Christoffel symbol  $\Gamma_{AB}^C$ is \emph{Darboux flat}
if there exist Darboux coordinates $z^A=(x^a,\theta_b)$ such that in these Darboux coordinates
the Christoffel symbol  $\Gamma_{AB}^C$ vanish, i.e., $\nabla_A\p_B=0$.
(A Darboux flat affine connection on $E$ induces the Darboux flat connection
$\g\colon \gamma_A=(-1)^B\Gamma^B_{AB}$ on volume forms.)

\begin{proposition}\label{prophalfdens}
The  half-density $\ss_{_C}(\nabla)$ does not depend on a choice of a connection in the class of the
Darboux flat connections: $\ss_{_C}( \nabla)=\ss_{_C}( \nabla')$
for two arbitrary Darboux flat affine connections $\nabla$ and $\nabla^\pr$.
\end{proposition}

This statement was used in a non-explicit way   in work~\cite{Kh2} where the half-density $\ss_{_C}(\nabla)$
was constructed in Darboux coordinates.

Proposition~\ref{prophalfdens} implies the existence of a canonical half-density on hypersurfaces in an
odd symplectic supermanifold.  This semi-density was first obtained by a straightforward calculation in~\cite{KhMkrt}.
On one hand, the invariant semi-density in an odd symplectic supermanifold is an analogue of Poincar\'{e}-Cartan integral invariants. On the other hand, the constructions above  are related with mean curvature of hypersurfaces (surfaces of codimension 1)
in the even Riemannian case. If $C$ is a surface of codimension $1|0$ in a Riemannian manifold $M$, then one can consider the canonical Levi-Civita connection  the and canonical volume form.  Applying constructions similar to the above we come to mean curvature of $C$.
In the odd symplectic case, there is no preferred affine connection compatible with the symplectic
structure (see \cite{Kh2} for details).

\appendix
\section{Connections and upper connections on volume forms}

{\small

A connection $\nabla$ on the bundle of volume forms defines a covariant derivative of volume forms with respect to vector fields. It obeys the natural linearity properties and the Leibniz rule:
\begin{itemize}
\item $\nabla_{\X}\left(\rh_1+\rh_2\right)=\nabla_{\X}\left(\rh_1\right)+\nabla_{\X}\left(\rh_2\right)$\,,
\item  $ \nabla_{f\X+g\Y}\left(\rh\right)=f\nabla_{\X}\left(\rh\right)+g\nabla_{\Y}\left(\rh\right)$\,,
\item $\nabla_{\X}\left(f \rh\right)=(\p_{\X}f)\,\rh +f\nabla_{\X}\left(\rh\right)$\,,
\end{itemize}
for   arbitrary volume forms $\rh$, $\rh_1$ and $\rh_2$, arbitrary vector fields $\X$ and $\Y$, and arbitrary functions $f$ and $g$.
Here $\p_\X$ is the ordinary derivative of a function  along a vector field.

\smallskip

Denote by $\nabla_a$ the covariant derivative with respect to the vector field $\p_a=\p/\p x^a$.   We have, for an arbitrary volume form $\rh=\rho(x) |Dx|$,
\begin{equation}
    \nabla_a\rh=\bigl(\p_a\rho +\gamma_a \rho\bigr)|D(x)|, \quad \text{where $\gamma_a(x) |D(x)|=\nabla_a(|D(x)|)$}\,.
\end{equation}
Under a change of local coordinates $x^a=x^a(x')$, the symbol
$\gamma_a$ transforms in the following way:
\begin{equation}\label{transformationofgenuineconnection}
   \gamma_a=x^{a'}_a \left(\gamma_{a'}+\p_{a'}\log \Bigl|\det \frac{\p x}{\p x'}\Bigr| \right)=
   x_a^{a'}\gamma_{a'}-x^b_{b'}x^{b'}_{ab}\,.
\end{equation}
We  use the shorthand notations for partial derivatives:
 $x^{a'}_a= \lder{x^{a'}\!\!}{x^a}$ and
 $x^{a'}_{bc}=\p^2 x^{a'}\!\!/\p x^b \p x^c$. Summation over repeated indices is assumed. Indeed,
\begin{multline*}
\gamma_a|D(x)|=\nabla_a |D(x)| = x^{a'}_a\nabla_{a^\pr}\left(\Bigl|\det\der{x}{x'}\Bigr|\cdot|D(x^\pr)|\right)=\\
x^{a'}_a\left(\p_{a'}\Bigl|\det\der{x}{x'}\Bigr|\cdot|D(x^\pr)|+ \Bigl|\det\der{x}{x'}\Bigr|\,\nabla_{a'}|D(x^\pr)|\right) =\\
x^{a'}_a\left(\p_{a'}\Bigl|\det\der{x}{x'}\Bigr|\cdot|D(x^\pr)|+ \gamma_{a'}\Bigl|\det\der{x}{x'}\Bigr|\,|D(x^\pr)|\right) =
x^{a'}_a\left(\p_{a'}\log\Bigl|\det\der{x}{x'}\Bigr| +
            \gamma_{a'}\right)\D{x}\,,
\end{multline*}
and we may also use  the standard relation $\delta\log \det M={\rm Tr\,}(M^{-1} \delta M)$.

Let  $S^{ab}$ be a contravariant tensor field. One can consider
a \emph{contravariant derivative} or  an \emph{upper connection}  $\uppernabla$  associated with $S$ (for an arbitrary vector bundle). This notion can be defined by axioms similar to those for a usual connection. We skip that and write down instead a coordinate expression for the particular case interesting for us. On volume forms, we   have
\begin{equation}\label{upperconnection1}
         \uppernabla^a\rh=\uppernabla^a\bigl(\rho(x) |D(x)|\bigr)=\left(S^{ab}\p_b\rho+\gamma^a\rho\right)|D(x)|\,.
\end{equation}

\begin{remark}
 Given a contravariant tensor field $S^{ab}$, a  connection $\nabla$ (covariant derivative) induces an upper connection (contravariant derivative)  $\uppernabla$ by the rule   $\uppernabla^a=S^{ab}\nabla_b$.
 If the  tensor field $S^{ab}$ is non-degenerate, the converse  is  also true.
 A non-degenerate contravariant tensor field $S^{ab}(x)$ induces a one-to-one correspondence between
upper connections and usual connections.
(Compare with the example \ref{operatorwothdegeneratesymbol}
below where an upper connection  does not in general define a connection.)
\end{remark}

Under a change of coordinates the symbol $\gamma^a$ for an upper connection \eqref{upperconnection1} transforms
as follows:
\begin{equation}\label{transformationofupperconnection}
      \gamma^{a'}=x^{a^\pr}_{a}\left(\gamma^a+S^{ab}\p_b\log \left|\det\der{x'}{x}\right|\right)\,.
\end{equation}
\begin{remark}
With some abuse of language we identify a connection (covariant derivative) on volume forms $\nabla$ with the symbol $\g=\{\gamma_a\}$: $\nabla_a\D{x}=\gamma_a\D{x}$.
\end{remark}

\begin{remark} It is worth noting  that the difference of two connections on volume forms is a
covector field, the difference of two upper connections on volume forms is
a vector field.  In other words the space of all connections (upper connections)
is an affine space associated with the linear space of the covector (vector) fields.
\end{remark}

Consider two important examples of connections on volume forms.

\begin{example}\label{connectioncorrespondingtovolumeform}
An arbitrary non-vanishing volume form $\rh_0$ (which can be regarded as a frame in the bundle of volume forms) defines a connection $\g^{\rh_0}$ by   the formula $\gamma_a=\gamma^{\rh_0}_a=-\p_a\log \rho_0(x)$.
  This is a flat connection: its curvature vanishes: $F_{ab}=\p_a\gamma_b-\p_a\gamma_a=0$.
  (Connection $\nabla$ considered in the formula \eqref{exampleofflatconnection}
is a flat connection defined by the volume form $\rh_g$.)
\end{example}

\begin{example}\label{secondexample}Let $\nabla^{TM}$ be an affine connection   on a
manifold $M$ (i.e., a connection on the tangent bundle). It defines a connection on volume forms $\nabla=-{\rm Tr\,} \nabla^{TM}$
with            $\gamma_a=-\Gamma^b_{ab}$
where $\Gamma^a_{bc}$ is the Christoffel symbol for $\nabla^{TM}$.
\end{example}

It is easy to see that a connection or an upper connection on volume forms define a covariant derivative  or, respectively, a contravariant
derivative on densities of  arbitrary weight: for a density $\ss=s(x)\D{x}^\l\in\F_\l$\,,
                      \begin{equation*}
        \nabla_a\ss=\bigl(\p_a s +\l\gamma_a s\bigr)\D{x}^\l\,.
                      \end{equation*}
Respectively, for an upper connection,
                     \begin{equation}\label{operationondensitiesofhigherweightsupper1}
     \nabla^a\ss=\left(S^{ab}\p_b s +\l\gamma^a s\right)\D{x}^\l\,.
                      \end{equation}
Sometimes we   use the concept of a \emph{connection of weight} $\delta$. This is a linear operation
that transforms  densities of weight $\l$  to  densities of weight $\mu=\l+\delta$:
for $\ss=s(x)\D{x}^\l\in\F_\l$\,,
\begin{equation*}
        \nabla_a\ss=\left(\p_a s +\l\gamma_a s\right)\D{x}^{\l+\delta}\,, \quad \text{where $\nabla_a\bigl(\D{x}\bigr)=\gamma_a\D{x}^{\l+\delta}$}\,.
\end{equation*}
Respectively, for an \emph{upper connection of weight } $\delta$,
\begin{equation*}
     \nabla^a\ss=\left(S^{ab}\p_b s  +\l\gamma^a s\right)\D{x}^{\l+\delta}\,, \quad \text{where $\nabla^a\bigl(\D{x}\bigr)=\gamma^a\D{x}^{\l+\delta}$}\,.	
\end{equation*}

}


\end{document}